\documentclass[12pt]{article}

\usepackage{amssymb}

\newtheorem{theorem}{Theorem}
\newtheorem{definition}[theorem]{Definition}
\newtheorem{lemma}[theorem]{Lemma}
\newtheorem{proposition}[theorem]{Proposition}
\newtheorem{corollary}[theorem]{Corollary}

\newtheorem{remark}[theorem]{\sc Remark}

\title{A survey on the inverse integrating factor.\thanks{The
authors are partially supported by a MCYT/FEDER grant number
MTM2008-00694}}

\author{{\sc Isaac A. Garc\' \i a$^{\ (1)}$ \& Maite Grau$^{\ (1)}$}}
\date{}

\begin{document}

\maketitle

\begin{abstract}
The relation between limit cycles of planar differential systems
and the inverse integrating factor was first shown in an article
of Giacomini, Llibre and Viano appeared in 1996. From that moment
on, many research articles are devoted to the study of the
properties of the inverse integrating factor and its relation with
limit cycles and their bifurcations. This paper is a summary of
all the results about this topic. We include a list of references
together with the corresponding related results aiming at being as
much exhaustive as possible. The paper is, nonetheless,
self-contained in such a way that all the main results on the
inverse integrating factor are stated and a complete overview of
the subject is given. Each section contains a different issue to
which the inverse integrating factor plays a role: the
integrability problem, relation with Lie symmetries, the center
problem, vanishing set of an inverse integrating factor,
bifurcation of limit cycles from either a period annulus or from a
monodromic $\omega$-limit set and some generalizations.
\end{abstract}
{\small{\noindent 2000 {\it AMS Subject Classification:} 34C07, 37G15, 34-02.   }}

\section{The Euler integrating factor \label{sect1}}

The method of integrating factors is, in principle, a means for
solving ordinary differential equations of first order and it is
theoretically important. The use of integrating factors goes back
to Leonhard Euler.

Let us consider a first order differential equation and write the
equation in the Pfaffian form
\begin{equation} \label{survey-1}
\omega = P(x,\,y)\,dy-Q(x,\,y)\,dx = 0 \ .
\end{equation}
We assume that the functions $P$ and $Q$ are of class
$\mathcal{C}^1$ in a region $\mathcal{U} \subseteq \mathbb{R}^2$.
If there is a solution of (\ref{survey-1}) which may be expressed
in the form $H(x,\,y) = h$ with $H$ having continuous partial
derivatives in $\mathcal{U}$ and with $h$ an arbitrary constant,
then it is not difficult to see that such an $H$  satisfies the
linear partial differential equation
\begin{equation} \label{survey-2}
P \frac{\partial H}{\partial x} + Q \frac{\partial H}{\partial y} = 0 \ .
\end{equation}
Conversely, every non-constant solution $H$ of (\ref{survey-2})
gives also a solution $H(x,\,y) = h$ of (\ref{survey-1}). Thus,
solving (\ref{survey-1}) and solving (\ref{survey-2}) are
equivalent tasks.

It is straightforward to show that if $H_0(x,\,y)$ is a
non-constant solution of equation (\ref{survey-2}), then all
solutions of this equation are of the form $F(H_0(x,\,y))$ where
$F$ is a freely chosen function with continuous derivative. The
connection between equations (\ref{survey-1}) and (\ref{survey-2})
may be presented also in another form. Suppose that $H(x,\,y) = h$
is any solution of (\ref{survey-1}). Then (\ref{survey-2}) implies
$$
\frac{\partial H / \partial y}{P} = -\frac{\partial H / \partial x}{Q} \ .
$$
If we denote the common value of these two ratios by $\mu(x,\,y)$,
then we have $\partial H / \partial y = \mu\, P$ and $\partial H /
\partial x = - \mu \, Q$. This gives to the differential of the
function $H$ the expression $d\,H(x,\,y) =
\mu(x,\,y)(P(x,\,y)\,dy-Q(x,\,y)\,dx)$. Hence, $\mu(x,\,y)$ is
called the integrating factor of the given differential equation
(\ref{survey-1}) because the left hand side of (\ref{survey-1})
turns, when multiplied by $\mu(x,\,y)$, to be an exact
differential.

Conversely, any integrating factor $\mu$ of (\ref{survey-1}), i.e.
such that $\mu(x,\,y)(P(x,\,y)\,dy-Q(x,\,y)\,dx)$ is the
differential of some function $H$, is easily seen to determine the
solutions of the form $H(x,\,y) = h$ of (\ref{survey-1}).
Altogether, solving the differential equation (\ref{survey-1}) is
equivalent to finding an integrating factor of the equation.

When an integrating factor $\mu$ of (\ref{survey-1}) is available,
the function $H$ can be obtained from the line integral
$$
H(x,\,y) = \int_{(x_0,\,y_0)}^{(x,\,y)} \mu(x,\,y)(P(x,\,y)\,dy-Q(x,\,y)\,dx)
$$
along any curve connecting an arbitrarily chosen point
$(x_0,\,y_0)$ and the point $(x,\,y)$ in the region $\mathcal{U}$.
We remark that this line integral might not be well-defined if the
region $\mathcal{U}$ is not simply-connected. When we know an
integrating factor $\mu$ of (\ref{survey-1}), we have a first
integral well-defined in each simply-connected subcomponent of the
region $\mathcal{U}$.

\section{The inverse integrating factor}

Let us consider a real planar autonomous differential system
\begin{equation}
\dot{x} \, = \, P(x,y), \qquad \dot{y} \, = \, Q(x,y), \label{eq1}
\end{equation}
where $P(x,y)$ and $Q(x,y)$ are of class
$\mathcal{C}^1(\mathcal{U})$ and $\mathcal{U} \subseteq
\mathbb{R}^2$ is an open set. The dot denotes derivation with
respect to the independent variable $t$ usually called {\em time},
that is $\dot{}=\frac{d}{dt}$. \par As usual, we associate to
system (\ref{eq1}) the vector field $\mathcal{X} \, = \, P(x,y)
\partial_x \, + \, Q(x,y) \partial_y$. Notice that the ordinary
differential equation $\omega = 0$ given in (\ref{survey-1}) is
just the differential equation of the orbits of system
(\ref{eq1}).

\begin{definition}
A function $V : \mathcal{U} \to \mathbb{R}$ is said to be an
inverse integrating factor of system {\rm (\ref{eq1})} if it is of
class $\mathcal{C}^1(\mathcal{U})$, it is not locally null and it
satisfies the following partial differential equation:
\begin{equation}\label{def-V}
P(x,y) \, \frac{\partial V(x,y)}{\partial x} \, + \,
Q(x,y) \, \frac{\partial V(x,y)}{\partial y} \, = \, \left(
\frac{\partial P(x,y)}{\partial x} \, + \, \frac{\partial
Q(x,y)}{\partial y} \right) \, V(x,y).
\end{equation}
\end{definition}

In short notation, an inverse integrating factor $V$ of system
(\ref{eq1}) satisfies $\mathcal{X} V = V {\rm div} \mathcal{X}$,
where ${\rm div} \mathcal{X} = \frac{\partial P}{\partial x} \, +
\, \frac{\partial Q}{\partial y}$ stands for the divergence of the
vector field $\mathcal{X}$.

Of course, the computation of an inverse integrating factor for a concrete system is a delicate
matter whose difficulty is comparable to solving the system itself.

If $V$ is an inverse integrating factor of a $\mathcal{C}^1$
vector field $\mathcal{X}$, then the zero set of $V$, $V^{-1}(0)
:=\{(x,y) \mid V(x,y)=0 \}$, is composed of trajectories of
$\mathcal{X}$. For by the equation (\ref{def-V}) that defines $V$,
$\mathcal{X}$ is orthogonal to the gradient vector field $\nabla
V$ along the zero set of $V$.

The name ``inverse integrating factor'' arises from the fact that
if $V$ solves equation (\ref{def-V}), then its reciprocal $1/V$ is
an integrating factor for $\mathcal{X}$ on $\mathcal{U} \setminus
V^{-1}(0)$.

\section{Local nontrivial Lie symmetries and inverse integrating factors \label{sect3}}

Roughly speaking, a symmetry group of a system of differential
equations is a continuous group which transforms solutions of the
system to other solutions. Simple typical examples are groups of
translations, rotations and scalings, but these certainly do not
exhaust the range of possibilities. Once one has determined the
symmetry group of a system of differential equations, a number of
applications become available.
\newline

More precisely, a symmetry of system (\ref{eq1})in $\mathcal{U}$,
where $\mathcal{U} \subseteq \mathbb{R}^2$ is an open set, is a
1--parameter Lie group of diffeomorphisms $\Phi_\epsilon$ acting
in $\mathcal{U}$ that maps the set of orbits of (\ref{eq1}) into
itself. When $\Phi_\epsilon(x,y) = (\bar{x}(x,y; \epsilon),
\bar{y}(x,y; \epsilon))$, the symmetry condition of (\ref{eq1})
reads for $\dot{\bar{x}} = P(\bar{x}, \bar{y})$, $\dot{\bar{y}} =
Q(\bar{x}, \bar{y})$ for all $\epsilon$ close to zero. Let the
$\mathcal{C}^1(\mathcal{U})$ vector field $\mathcal{Y} = \xi(x, y)
\partial_x + \eta(x, y) \partial_y$ be the infinitesimal generator
of the 1--parameter Lie group $\Phi_\epsilon$, that is,
$\bar{x}(x,y;\epsilon) =  x+ \epsilon \xi(x,y) + O(\epsilon^2)$,
$\bar{y}(x,y;\epsilon) = y+ \epsilon \eta(x,y) + O(\epsilon^2)$.
Denoting by $\mathcal{X} \, = \, P(x,y) \partial_x \, + \, Q(x,y)
\partial_y$ the vector field associated to system (\ref{eq1}), it
is well known that a characterization of the Lie symmetries of
(\ref{eq1}) is given by the relation $[\mathcal{X}, \mathcal{Y}] =
\mu(x,y) \mathcal{X}$ for certain scalar function $\mu :
\mathcal{U} \to \mathbb{R}$. In this expression we have used the
{\it Lie bracket} of two $\mathcal{C}^1$-vector fields ${\cal X}$
and ${\cal Y}$ defined as $[{\cal X}, {\cal Y} ]:= {\cal X} {\cal
Y} - {\cal Y} {\cal X}$. Using coordinates we have {\small
\begin{equation}
[{\cal X}, {\cal Y} ] = \left( P {\partial \xi \over \partial x}- \xi {\partial P \over \partial x}+ Q {\partial \xi \over \partial y} - \eta {\partial P \over \partial y} \right) \partial_x + \left( P {\partial \eta \over \partial x}- \xi {\partial Q \over \partial x}+ Q {\partial \eta \over \partial y} - \eta {\partial Q \over \partial y} \right) \partial_y \ .
\label{sim2.4.1}
\end{equation}
}
When beginning students first encounter ordinary differential equations, they are presented with a variety of special techniques designed to solve certain particular types of equations, such as separable, homogeneous or exact. Indeed, this was the state of the art around the middle of the nineteenth century, when Sofus Lie made the profound discovery that these special methods were, in fact, all special cases of a general integration procedure based on the invariance of the differential equation under a continuous group of symmetries. This observation at once unified and significantly extended the available integration techniques.

\begin{center}
\begin{tabular}{||c|c||}
\hline\hline \\
{\sc Differential Equation}  & {\sc Lie Symmetry }
\\ \hline\hline $dy/dx = f(x)g(y)$ & ${\cal Y} = g(y) \partial_y$ \\
\hline $dy/dx = f(a x+b y)$ & ${\cal Y} = b \partial_x + a \partial_y$ \\ \hline $dy/dx =
\frac{y+x f(\sqrt{x^2+y^2})}{x-y f(\sqrt{x^2+y^2})}$ & ${\cal Y} =
y \partial_x - x \partial_y$ \\ \hline $dy/dx =
f(y/x)$ & ${\cal Y} = x \partial_x + y \partial_y$ \\ \hline
$dy/dx =  P(x) y+Q(x)$ & ${\cal Y} = \exp\left( \int P(x) dx \right)
\partial_y$ \\ \hline
$dy/dx =  P(x) y+Q(x) y^n$ & ${\cal Y} = y^n \exp\left[ (1-n) \int
P(x) dx \right] \partial_y$ \\ \hline
\end{tabular}
\end{center}

Consider now a $\mathcal{C}^1$ vector field ${\cal X} = P(x,y)
\partial_x + Q(x,y) \partial_y$ defined in an open connected
subset $\mathcal{U}  \subseteq \mathbb{R}^2$. In the case of a
single first order ordinary differential equation $d y / d x =
Q(x,y) / P(x,y)$, the Lie symmetries method provides by quadrature
an explicit formula for the general solution. In fact, one can
easily see that if we know a Lie symmetry in $\mathcal{U}$ with
infinitesimal generator ${\cal Y} = \xi(x,y) \partial_x +
\eta(x,y) \partial_y$ then we construct an inverse integrating
factor $V = \det \{{\cal X}, {\cal Y}\} = P \eta - Q \xi$ defined
in $\mathcal{U}$, but the converse is not always true. To see that, assume now
the existence of an inverse integrating factor $V$ of $\mathcal{X}$ in a simply connected domain $\mathcal{U}$ and we look for an  infinitesimal generator $\mathcal{Y} = \xi(x,y) \partial_x + \eta(x,y)
\partial_y$ of a Lie symmetry  of $\mathcal{X}$ well defined in
$\mathcal{U}$. We recall that a singular point $p \in U$ of $\mathcal{X}$ is
called {\it weak} if ${{\textrm{div}}} \mathcal{X}(p) = 0$. If 
there is no weak singularity of $\mathcal{X}$ in $\mathcal{U}$, then we can
do at least one of the following constructions:
\begin{description}
\item[(i)] Prescribe the function $\xi(x,y)$ and solve $\eta(x,y)$
from $V= P \eta- Q \xi$.

\item[(ii)] Prescribe the function $\eta(x,y)$ and solve
$\xi(x,y)$ from $V= P \eta- Q \xi$.

\item[(iii)] Take the rescaled hamiltonian vector field
$$
\mathcal{Y}= \frac{1}{{{\textrm{div}}} \mathcal{X}} \
(-\frac{\partial{V}}{\partial{y}}\partial_x +
\frac{\partial{V}}{\partial{x}}\partial_y) \ ,
$$
defined in $\mathcal{U}  \backslash \{ (x,y) \in  \mathcal{U} : {{\textrm{div}}}
\mathcal{X} = 0 \}$.
\end{description}
Therefore, the equivalence between inverse integrating factors and
Lie symmetries for planar vector fields $\mathcal{X}$ is not true,
in general, in neighborhoods of weak singular points of
$\mathcal{X}$. Of course, some special situations can appear
giving the equivalence when $\mathcal{X}$ possesses an analytic
first integral in these neighborhoods as the nondegenerate center
singular point shows.
\newline

Importance of inverse integrating factors arises
from the fact that the differential 1-form $\omega / V = (P \ dy -
Q \ dx) / V$ is closed ($d (\omega/V) =0$) in $\mathcal{U}
\backslash V^{-1}(0)$. Then in the case in which $\mathcal{U}
\backslash V^{-1}(0)$ is simply-connected, the 1-form $\omega / V$
is {\it exact} ($\omega/V = d H$), and therefore a ${\cal C}^2$
first integral $H(x,y)$ of the differential equation is
immediately constructed. As a consequence, the vector field
$\mathcal{X} \, = \, P(x,y)
\partial_x \, + \, Q(x,y) \partial_y$ is topologically equivalent,
in $\mathcal{U}$, to the hamiltonian vector field $\mathcal{X} / V
\, = \, \frac{\partial H}{\partial y}
\partial_x \, - \, \frac{\partial H}{\partial x} \partial_y$.
\newline

Making a pause in this exposition we now present an example. Let us consider the following cubic system
\begin{equation}
\dot{x}= P(x,y) = -y-x(x^2+y^2-1) \ , \ \ \dot{y}= Q(x,y)= x-y(x^2+y^2-1) \ .
\label{ej}
\end{equation}
An inverse integrating factor for system (\ref{ej}) is given by $V(x,y)=(x^2+y^2)(x^2+y^2-1)$. Associated to him one has the first integral
$$
H(x,y)={(x^2+y^2-1) \over (x^2+y^2)} \, \exp\left\{ 2 \, \arctan
\left({y \over x}\right) \right\}  ,
$$
which is not continuous in $(0,0)$. On the other hand, since the
polar form of the system is $\dot{r} = 2 r^2 (r^2-1)$,
$\dot{\varphi} = 1$ it is easy to check that the unit circle
$x^2+y^2-1=0$ is the unique limit cycle of system (\ref{ej}). Let
${\cal X}=P(x,y) \partial / \partial x + Q(x,y) \partial /
\partial y$ be the vector field associated with system (\ref{ej}).
From the symmetries point of view, since ${\cal Y}=y  \partial_x
-x \partial_y$ satisfies $[ {\cal X} , {\cal Y} ] \equiv 0$ we
have that ${\cal Y}$ is the infinitesimal generator of a Lie group
admitted by system (\ref{ej}) which is just the $SO(2)$ rotation
group $\bar{x} = x \cos\epsilon - y \sin\epsilon$, $\bar{y} = x
\sin\epsilon + y \cos\epsilon$. Hence $V(x,y) = \det \{{\cal X},
{\cal Y}\}$ is an inverse integrating factor of system (\ref{ej}).
Notice that the only common integral curves for the vector fields
${\cal X}$ and ${\cal Y}$ are included in $V^{-1}(0)$ and are just
the separatrices of ${\cal X}$. This behavior will be explained in
future sections.

\section{On the integrability problem \label{sectinteg}}

The integrability problem is mainly related to planar {\bf
polynomial} differential systems of the form
\begin{equation}
\dot{x}=P(x,y), \quad \dot{y}=Q(x,y), \label{0eq1}
\end{equation}
where $P(x,y), Q(x,y) \in \mathbb{R}[x,y]$ are coprime
polynomials, that is, there is no non-constant polynomial which
divides both $P$ and $Q$. We call {\rm d} the maximum degree of
$P$ and $Q$ and we say that system (\ref{0eq1}) is of degree {\rm
d}. When ${\rm d}=2$, we say that (\ref{0eq1}) is a {\em quadratic
system}. \par If $p$ is a point such that $P(p)=Q(p)=0$, then we
say that $p$ is a {\em singular point} of system (\ref{0eq1}).
\newline

As we have already defined in Section \ref{sect1}, a
${\mathcal{C}}^j$ function $H: {\mathcal{U}} \to \mathbb{R}$ such
that it is constant on each trajectory of (\ref{0eq1}) and it is
not locally constant is called a {\em first integral} of system
(\ref{0eq1}) of class $j$ defined on $\mathcal{U} \subseteq
\mathbb{R}^2$. The equation $H(x,y)=h$ for a fixed $h \in
\mathbb{R}$ gives a set of trajectories of the system, but in an
implicit way. When $j \geq 1$, these conditions are equivalent to
$P(x,y) \frac{\partial H}{\partial x} + Q(x,y) \frac{\partial
H}{\partial y} = 0$ and $H$ not locally constant. The problem of
finding such a first integral and the functional class it must
belong to is what we call the {\em integrability problem}.
\par
To find an integrating factor or an inverse integrating factor for
system (\ref{0eq1}) is closely related to finding a first integral
for it. When considering the integrability problem we are also
addressed to study whether an (inverse) integrating factor belongs
to a certain given class of functions.
\newline

When a first integral $H$ of system (\ref{0eq1}) is known, all the
orbits of the system are contained in its domain of definition are
given by the level sets $H(x,y)=h$. Thus, a natural strategy is to
look for the determination of some of the orbits of the system and
try to build a first integral with them. In particular, and since
system (\ref{0eq1}) is polynomial, those orbits which are
algebraic will be of special interest. \par An {\em invariant
curve} is a curve given by $f(x,y)=0$, where $f: \mathcal{U}
\subseteq \mathbb{R}^2 \to \mathbb{R}$ is a $\mathcal{C}^1$
function in the open set $\mathcal{U}$, non locally constant and
such that there exists a $\mathcal{C}^1$ function in
$\mathcal{U}$, denoted by $k(x,y)$ and called {\em cofactor},
which satisfies:
\begin{equation}
P(x,y) \, \frac{\partial f}{\partial x} (x,y) \, + \, Q(x,y) \,
\frac{\partial f}{\partial y}(x,y) \, = \, k(x,y) \, f(x,y),
\label{0ci}
\end{equation}
for all $(x,y) \in \mathcal{U}$. The notion of invariant curve was
first introduced in \cite{GarciaGine}. The identity (\ref{0ci})
can be rewritten by $\mathcal{X} f = k f$. We recall that
$\mathcal{X} f$ denotes the scalar product of the vector field
$\mathcal{X}$ and the gradient vector $\nabla f$ related to
$f(x,y)$, that is, $\nabla f(x,y)=( \frac{\partial f}{\partial x}
(x,y) , \frac{\partial f}{\partial y} (x,y) )$. We will denote by
$\frac{d f}{d t}$ or by $\dot{f}$ the function $\mathcal{X} f$
once evaluated on a solution of system (\ref{0eq1}). In case
$f(x,y)=0$ defines a curve in the real plane, this definition
implies that the function $\mathcal{X} f$ is equal to zero on the
points such that $f(x,y)=0$. In the article \cite{GarciaGine} an
invariant curve is defined as a $\mathcal{C}^1$ function $f(x,y)$
defined in the open set $\mathcal{U} \subseteq \mathbb{R}^2$, such
that, the function $\mathcal{X} f$ is zero in all the points $\{
(x,y) \in \mathcal{U} \, | \, f(x,y)=0 \}$. We notice that our
definition of invariant curve is a particular case of the previous
one but, for the sake of our results, the cofactor is very
important and that's why we always assume its existence.
\par
When the cofactor $k(x,y)$ is a polynomial, we say that $f(x,y)=0$
is an invariant curve with polynomial cofactor. We only admit
invariant curves with polynomial cofactor of degree lower or equal
than ${\rm d}-1$, that is $\deg k(x,y) \leq {\rm d}-1$, where
${\rm d}$ is the degree of system (\ref{0eq1}).
\par
The notion of invariant curve is a generalization of the notion of
invariant algebraic curve.  An {\em invariant algebraic curve} is
an algebraic curve $f(x,y)=0$, where $f(x,y) \in \mathbb{C}[x,y]$,
which is invariant by the flow of system (\ref{0eq1}). This
condition equals to $\mathcal{X} f = k f$, where the cofactor of
an invariant algebraic curve is always a polynomial of degree
$\deg k(x,y) \leq {\rm d}-1$.
\par
We cite \cite{Llibre, Schlomiuk1, Schlomiuk2} as compendiums of
the results on invariant algebraic curves. For instance, in
\cite{Llibre}, it is shown that if $f(x,y)=0$ and $g(x,y)=0$ are
two invariant algebraic curves of system (\ref{0eq1}) with
cofactors $k_f(x,y)$ and $k_g(x,y)$, respectively, then the
product of the two polynomials gives rise to the curve $(f g)(x,y)
=0$ which is also an invariant algebraic curve of system
(\ref{0eq1}) and whose cofactor is $k_f(x,y) + k_g(x,y)$.
\par
In order to state the known results of integrability using
invariant algebraic curves, we need to consider complex algebraic
curves $f(x,y)=0$, where $f(x,y) \in \mathbb{C}[x,y]$. Since
system (\ref{0eq1}) is defined by real polynomials, if $f(x,y)=0$
is an invariant algebraic curve with cofactor $k(x,y)$, then its
conjugate $\bar{f}(x,y)=0$ is also an invariant algebraic curve
with cofactor $\bar{k}(x,y)$. Hence, its product $f(x,y)
\bar{f}(x,y) \in \mathbb{R}[x,y]$ gives rise to a real invariant
algebraic curve with a real cofactor $k(x,y)+ \bar{k}(x,y)$. For a
sake of simplicity, we consider invariant algebraic curves defined
by polynomials in $\mathbb{C}[x,y]$, although we always keep in
mind the previous observation. In $\mathbb{R}^2$, the curve given
by $f(x,y)=0$, where $f(x,y)$ is a real function, may only contain
a finite number of isolated singular points or be the null set.
\par
An algebraic curve $f(x,y)=0$ is called {\em irreducible} when
$f(x,y)$ is an irreducible polynomial in the ring
$\mathbb{C}[x,y]$. We can assume, without loss of generality, that
$f(x,y)$ is an irreducible polynomial in $\mathbb{C}[x,y]$,
because if $f(x,y)$ is reducible, then all its proper factors give
rise to invariant algebraic curves. Given an algebraic curve
$f(x,y)=0$, we can always assume that the polynomial $f(x,y)$ has
no multiple factors, that is, its decomposition in the ring
$\mathbb{C}[x,y]$ is of the form $f(x,y) = f_1(x,y) f_2(x,y)
\ldots f_{\ell}(x,y)$, where $f_i(x,y)$ are irreducible
polynomials and $f_i(x,y) \neq c f_j(x,y)$ if $i \neq j$ and for
any $c \in \mathbb{C}$. The assumption that given an algebraic
curve $f(x,y)=0$, the polynomial $f(x,y)$ has no multiple factors
is mainly used to ensure that we do not consider ``false''
singular points. If $p$ is a point such that $f(p)=0$ and $\nabla
f(p)=0$, and $f(x,y)$ has no multiple factors, then $p$ is a
singular point of the curve $f(x,y)=0$. But, if $f(x,y)$ has
multiple factors, for instance, $f(x,y)=f_1(x,y)^2$ where
$f_1(x,y)$ is an irreducible polynomial in $\mathbb{C}[x,y]$, then
all the points of the curve $\{ p \, | \, \,  f_1(p)=0 \}$ satisfy
the property that $f(p)=0$ and $\nabla f(p) =0$ although they are
not all singular points.
\par
We recall that if $p$ is a singular point of an invariant
algebraic curve $f(x,y)=0$ of a system (\ref{0eq1}), then $p$ is a
singular point of the system. Given an algebraic curve $f(x,y)=0$,
we will always assume that the decomposition of $f(x,y)$ in the
ring $\mathbb{C}[x,y]$ has no multiple factors. We want to
generalize this property to invariant curves, that's why we will
always assume that, given an invariant curve $f(x,y)=0$, if $p \in
\mathcal{U}$ is such that $f(p)=0$ and $\nabla f(p)=0$, then $p$
is a singular point of system (\ref{0eq1}). This technical
hypothesis generalizes the notion of not having multiple factors
for algebraic curves. In \cite{seiden4}, a set of necessary
conditions for a system (\ref{0eq1}) to have an irreducible
invariant algebraic curve is given.
\newline

Invariant algebraic curves are the main objects used in the
Darboux theory of integrability. In \cite{Darboux}, G. Darboux
gives a method for finding an explicit first integral for a system
(\ref{0eq1}) in case that ${\rm d}({\rm d}+1)/2 + 1$ different
irreducible invariant algebraic curves are known, where {\rm d} is
the degree of the system. In this case, a first integral of the
form $H=f_1^{\lambda_1}f_2^{\lambda_2} \ldots f_s^{\lambda_s}, $
where each $f_i(x,y)=0$ is an invariant algebraic curve for system
(\ref{0eq1}) and $\lambda_i \in \mathbb{C}$ not all of them null,
for $i=1,2,\ldots,s$, $s \in \mathbb{N}$, can be constructed. The
functions of this type are called {\em Darboux functions}.
\par
As we have already stated, given an invariant algebraic curve
$f(x,y)=0$ whose imaginary part is not null, then its conjugate is
also an invariant algebraic curve. Moreover, as system
(\ref{0eq1}) is real, if $f(x,y)$ appears in the expression of a
first integral of the form given by Darboux with exponent
$\lambda$, then $\bar{f}(x,y)$ appears in the same expression with
exponent $\bar{\lambda}$. We call ${\rm Re} f$ the real part of
the polynomial $f$ and by ${\rm Im} f$ its imaginary part.
Analogously, let us call ${\rm Re} \lambda$ the real part of the
complex number $\lambda$ and by ${\rm Im} \lambda$ its imaginary
part. We call $\mathbf{i}=\sqrt{-1}$ and we use the following
formula for complex numbers:
\[ \arctan(z) = \log \left[ \left( \frac{1 - \mathbf{i} z}{1 + \mathbf{i} z} \right)
^{ \mathbf{i}/2} \right], \quad z \in \mathbb{C}, \] to show that
\begin{eqnarray*}
f^{\lambda} \bar{f}^{\bar{\lambda}} & = & \left( {\rm Re} f + {\rm
Im} f \, \mathbf{i} \right) ^{{\rm Re} \lambda + {\rm Im} \lambda
\, \mathbf{i}} \ \left( {\rm Re} f - {\rm Im} f\, \mathbf{i}
\right) ^{{\rm Re} \lambda - {\rm Im} \lambda \, \mathbf{i}}
\\ & = & \left( ({\rm Re} f)^2 + ({\rm Im} f)^2 \right)^{{\rm Re} \lambda} \ \exp \left\{
-2\,  {\rm Im} \lambda \, \arctan \left( \frac{{\rm Im} f}{{\rm
Re} f} \right) \right\}.
\end{eqnarray*}
We deduce that the product $f(x,y)^{\lambda}
\bar{f}(x,y)^{\bar{\lambda}}$ is a real function and so it is any
Darboux function $H=f_1^{\lambda_1}f_2^{\lambda_2} \ldots
f_s^{\lambda_s}$.
\par
We have that the Darboux function $H$ can be defined in the open
set $\mathbb{R}^2 \setminus \Sigma$, where $ \Sigma = \{ (x,y) \in
\mathbb{R}^2 \mid (f_1 \cdot f_2 \cdot \dots \cdot f_r)(x,y) = 0
\}.$ We remark that, particularly, if $\lambda_i \in \mathbb{Z}$ ,
$\forall i=1,2,\ldots,r$, $H$ is a {\em rational first integral}
for system (\ref{0eq1}). In this sense J. P. Jouanoulou
\cite{Jouanolou}, showed that if at least ${\rm d}({\rm d}+1) + 2$
different irreducible invariant algebraic curves are known, then
there exists a rational first integral.
\par
The main fact used to prove Darboux's theorem (and Jouanoulou's
improvement) is that the cofactor corresponding to each invariant
algebraic curve is a polynomial of degree $\leq {\rm d}-1$.
Invariant curves with polynomial cofactor can also be used in
order to find a first integral for the system. This observation
enables a generalization of the Darboux's theory which is given in
\cite{GarciaGine1}, where, for instance, non-algebraic invariant
curves with an algebraic cofactor for a polynomial system of
degree $4$ are presented. In \cite{rocky}, other examples are
given of such invariant curves with polynomial cofactor for some
families of systems and the way they are used to construct
explicit first integrals and inverse integrating factors for the
corresponding systems. As a continuation of \cite{rocky}, in
\cite{rocky2} we study when a planar differential system
polynomial in one variable linearizes in the sense that it has an
inverse integrating factor which can be constructed by means of
the solutions of linear differential equations and we describe
some families of differential systems which are Darboux integrable
and whose inverse integrating factor is constructed using the
solutions of a second--order linear differential equation defining
a family of orthogonal polynomials.\par Some generalizations of
the classical Darboux theory of integrability may be found in the
literature. For instance, independent singular points can be taken
into account to reduce the number of invariant algebraic curves
necessary to ensure the Darboux integrability of the system, see
\cite{ChLlSoto}. A good summary of many of these generalizations
can be found in \cite{Chara} and a survey on the integrability of
two-dimensional systems can be found in \cite{flows}. One of the
most important definitions in this sense is the notion of
exponential factor which is given by C. Christopher in
\cite{Christopher1}, when he studies the multiplicity of an
invariant algebraic curve. The notion of exponential factor is a
particular case of invariant curve for system (\ref{0eq1}). Given
two coprime polynomials $h,g \in \mathbb{R}[x,y]$, the function
$e^{h/g}$ is called an {\em exponential factor} for system
(\ref{0eq1}) if for some polynomial $k$ of degree at most ${\rm
d}-1$, where {\rm d} is the degree of the system, the following
relation is fulfilled:
\[ P \left( \frac{\partial\, e^{h/g}}{\partial x} \right) + Q
\left( \frac{\partial\, e^{h/g}}{\partial y} \right) = k(x,y) \,
\, e^{h/g}.
\] As before, we say that $k(x,y)$ is the {\em cofactor} of the
exponential factor $e^{h/g}$.
\newline

The next proposition, proved in \cite{Christopher1}, gives the
relationship between the notion of invariant algebraic curve and
exponential factor.
\begin{proposition} {\sc \cite{Christopher1}} \  If $F=e^{h/g}$ is
an exponential factor and $g$ is not a constant, then $g=0$ is an
invariant algebraic curve, and $h$ satisfies the equation $P
\frac{\partial h}{\partial x} + Q \frac{\partial h}{\partial y}= h
\, k_g + g\, k_F$ where $k_g$ and $k_F$ are the cofactors of $g$
and $F$, respectively.
\end{proposition}
The notion of exponential factor is very important in the Darboux
theory of integrability since it does not only allow the
construction of first integrals following the same method
described by Darboux, but it also explains the meaning of the
multiplicity of an invariant algebraic curve in relation with the
differential system (\ref{0eq1}). A complete work on this subject
can be found in \cite{ChLlPe}.
\par
In the same way as with invariant algebraic curves, given an
exponential factor $F=\exp\{ h/g \}$, since system (\ref{0eq1}) is
a real system, there is no lack of generality in considering that
$h(x,y), g(x,y) \in \mathbb{R}[x,y]$. If $F=\exp\{ h/g \}$ is an
exponential factor with non-null imaginary part, then its complex
conjugate, $\bar{F} = \exp \{ \bar{h} / \bar{g} \}$ is also an
exponential factor, as it can be easily checked by its defining
equation. Moreover, the product $F \, \bar{F}= \exp \{ h/g +
\bar{h}/\bar{g} \}$ is a real exponential factor with a real
cofactor.
\newline

Since the notion of exponential factor is the most current
generalization in the Darboux theory of integrability, any
function of the form:
\begin{equation}
f_1^{\lambda_1} f_2^{\lambda_2} \cdots f_r^{\lambda_r} \left( \exp
\left( \frac{h_1}{g_1^{n_1}} \right)\right)^{\mu_1} \left( \exp
\left( \frac{h_2}{g_2^{n_2}} \right)\right)^{\mu_2} \cdots \left(
\exp \left( \frac{h_\ell}{g_\ell^{n_\ell}}
\right)\right)^{\mu_\ell}, \label{0gdf}
\end{equation}
where $r, \ell \in \mathbb{N}$, $f_i(x,y)=0$ ($1 \leq i \leq r$)
and $g_j(x,y)=0$ ($1 \leq j \leq \ell$) are invariant algebraic
curves of system (\ref{0eq1}), $h_j(x,y)$ ($1 \leq j \leq \ell$)
are polynomials in $\mathbb{C}[x,y]$, $\lambda_i$ ($1 \leq i \leq
r$) and $\mu_j$ ($1 \leq j \leq \ell$) are complex numbers and
$n_j$ ($1 \leq j \leq \ell$) are non-negative integers, is called
a {\em (generalized) Darboux function}. 
\newline

Let us present a short survey about the Darboux method and its improvements. Let us recall
that a singular point $(x_0, y_0)$ of system (\ref{0eq1}) is
called {\em weak} if the divergence, ${\rm div}\mathcal{X}$, of
system (\ref{0eq1}) at $(x_0, y_0)$ is zero. We recall that
$\mathcal{X}$ denotes the vector field associated to system
(\ref{0eq1}). We denote by $\mathbb{C}_{\rm d-1} [x,y]$ the set of
polynomials in $\mathbb{C}[x,y]$ of degree lower than ${\rm d}$.
We say that $s$ points $(x_k, y_k)\in \mathbb{C}^2$,
$k=1,2,\ldots,s$, are {\em independent} with respect to
$\mathbb{C}_{\rm d-1} [x,y]$ if the intersection of the $s$
hyperplanes
\[ \left\{ \left( a_{ij}\right) \in \mathbb{C}^{\rm d(d+1)/2} \, :
\, \sum_{i+j=0}^{\rm d-1} x_k^i y_k^i a_{ij} \, = \,
0\right\}_{k=1,2,\ldots,s} \] is a linear subspace of
$\mathbb{C}^{\rm d(d+1)/2}$ of dimension ${\rm d(d+1)/2} -s > 0$.

The main results about the Darboux method and its improvements are summarized in the following theorem, which can be
found in \cite{LlibPan04}, see also \cite{Chara}.

\begin{theorem} \label{thLlibPan04}
Suppose that a polynomial differential system {\rm (\ref{0eq1})}
of degree ${\rm d}$ admits $r$ irreducible invariant algebraic
curves $f_i = 0$ with cofactors $K_i$ for $i = 1,2,\ldots,r$;
$\ell$ exponential factors $exp(h_j/g_j^{n_j})$ with cofactors
$L_j$ for $j =1,2,\ldots, \ell$; and $s$ independent singular
points $(x_k, y_k)$ such that $f_i(x_k, y_k) \neq 0$ for $i =
1,2,\ldots,r$ and for $k=1,2,\ldots,s$. Moreover, the irreducible
factors of the polynomials $g_j$ are some $f_i$'s.
\begin{itemize}
\item[{\rm (a)}] There exist $\lambda_i, \mu_j\in \mathbb{C}$ not all zero such that
$\sum_{i=1}^{r} \lambda_i K_i + \sum_{j=1}^{\ell} \mu_j L_j=0$, if
and only if the (multi–valued) function {\rm (\ref{0gdf})} is a
first integral of system {\rm (\ref{0eq1})}.
\item[{\rm (b)}] If $r + \ell + s = [{\rm d}({\rm d} + 1)/2] + 1$,
then there exist $\lambda_i, \mu_j\in \mathbb{C}$ not all zero
such that $\sum_{i=1}^{r} \lambda_i K_i + \sum_{j=1}^{\ell} \mu_j
L_j=0$.
\item[{\rm (c)}] If $r + \ell + s \geq [{\rm d}({\rm d} + 1)/2] + 2$,
then system {\rm (\ref{0eq1})} has a rational first integral, and
consequently all trajectories of the system are contained in
invariant algebraic curves.
\item[{\rm (d)}] There exist $\lambda_i, \mu_j\in \mathbb{C}$ not all zero such that
$\sum_{i=1}^{r} \lambda_i K_i + \sum_{j=1}^{\ell} \mu_j L_j={\rm
div}\, \mathcal{X}$ if and only if the function {\rm (\ref{0gdf})}
is an inverse integrating factor of system {\rm (\ref{0eq1})}.
\item[{\rm (e)}] If $r + \ell + s = {\rm d}({\rm d} + 1)/2$ and
$s$ independent singular points are weak, then the function {\rm
(\ref{0gdf})} for convenient $\lambda_i, \mu_j\in \mathbb{C}$ not
all zero is a first integral or an inverse integrating factor of
system {\rm (\ref{0eq1})}.
\end{itemize}
\end{theorem}

Introducing the notion of multiplicity of the invariant algebraic
hypersurfaces of a polynomial vector field in $\mathbb{C}^n$, the
results of Darboux integrability theory of Theorem
\ref{thLlibPan04} have been generalized to systems in
$\mathbb{C}^n$, where $n \geq 2$, see \cite{Lli-Zha} and the
references therein.
\newline

An improvement of the previous Darboux theorem is presented in
\cite{ChaGiaGin99} when the system has a center. As usual $\lfloor
q \rfloor$ means the integer part of the real number $q$.
\begin{theorem} {\rm \cite{ChaGiaGin99}} \label{thChaGiaGin99}
Consider a polynomial system {\rm (\ref{0eq1})} of degree ${\rm
d}$, with a center at the origin and with an arbitrary linear
part. Suppose that this system admits ${\rm d}({\rm d} + 1)/2 -
\lfloor ({\rm d} + 1)/2\rfloor$ invariant algebraic curves or
exponential factors. Then this system has a Darboux inverse
integrating factor.
\end{theorem}
In the following section we present several relations between the
existence of an inverse integrating factor and the center problem.
\newline

We recall that the integrability problem consists in finding the
class of functions a first integral of a given system (\ref{0eq1})
must belong to. We have system (\ref{0eq1}) defined in a certain
class of functions, in this case, the polynomials with real
coefficients $\mathbb{R}[x,y]$, and we consider the problem
whether there is a first integral in another, possibly larger,
class. For instance in \cite{Poin97}, H. Poincar\'e stated the
problem of determining when a system (\ref{0eq1}) has a rational
first integral. The works of M.J. Prelle and M.F. Singer
\cite{PrelleSinger} and M.F. Singer \cite{Singer} go on this
direction since they give a characterization of when a polynomial
system (\ref{0eq1}) has an elementary or a Liouvillian first
integral. An important fact of their results is that invariant
algebraic curves play a distinguished role in this
characterization. Moreover, this characterization is expressed in
terms of the inverse integrating factor.
\par
Roughly speaking, an {\em elementary function} is a function
constructed from rational functions by using algebraic operations,
composition and exponentials, applied a finite number of times,
and a {\em Liouvillian function} is a function constructed from
rational functions by using algebraic operations, composition,
exponentials and integration, applied a finite number of times. A
precise definition of these classes of functions is given in
\cite{PrelleSinger, Singer}. We are mainly concerned with
Liouvillian functions but we will state some results related to
integration of a system (\ref{0eq1}) by means of elementary
functions. \par We recall that $\mathbb{C}(x,y)$ denotes the
quotient field associated to the ring of polynomials with complex
coefficients, that is, $\mathbb{C}(x,y)$ is the field of rational
functions with complex coefficients.

\begin{theorem} {\sc \cite{PrelleSinger}} \ If the system {\rm (\ref{0eq1})} has
an elementary first integral, then there exist $w_0, w_1, \ldots,
w_n$ algebraic over the field $\mathbb{C}(x,y)$ and $c_1, c_2,
\ldots, c_n$ in $\mathbb{C}$ such that the elementary function
\begin{equation} H = w_0 + \sum_{i=1}^{n} c_i \ln (w_i) \label{0ef} \end{equation} is a first
integral of system {\rm (\ref{0eq1})}. \label{0thprellesinger0}
\end{theorem}
The existence of an elementary first integral is intimately
related to the existence of an algebraic inverse integrating
factor, as the following result shows.
\begin{theorem} {\sc \cite{PrelleSinger}} \ If the
system {\rm (\ref{0eq1})} has an elementary first integral, then
there is an inverse integrating factor of the form \[ V= \left(
\frac{A(x,y)}{B(x,y)} \right)^{1/N} , \] where $A, B \in
\mathbb{C}[x,y]$ and $N$  is an integer number.
\label{0thprellesinger}
\end{theorem}

The paper \cite{CGGLl2} is devoted to study which is the form of
the inverse integrating factor of a polynomial planar system
(\ref{0eq1}) with a Darboux first integral $H$ of the form
(\ref{0gdf}). This work is an improvement of the results of Prelle
and Singer in \cite{PrelleSinger} where it is shown that these
Darboux integrable vector fields have a rational inverse
integrating factor (see Theorem 7 of \cite{PrelleSinger}). In
\cite{CGGLl2}, another proof of this result is presented.
\begin{theorem} {\sc \cite{CGGLl2}} \ If the system {\rm
(\ref{0eq1})} has a (generalized) Darboux first integral of the
form {\rm (\ref{0gdf})}, then there is a rational inverse
integrating factor, that is, an inverse integrating factor of the
form:  \[ V= \frac{A(x,y)}{B(x,y)} , \] where $A, B \in
\mathbb{C}[x,y]$. \label{0thcggll2}
\end{theorem}
Unfortunately, not all the elementary functions of the form
(\ref{0ef}) are of (generalized) Darboux type. That's why, we can
find systems with an elementary first integral and without a
rational inverse integrating factor. The following example is of
this type. The system appears in the works of Jean
Moulin-Ollagnier \cite{Moulin1, Moulin}, although he does not give
an explicit expression for the first integral. The Lotka-Volterra
system:
\begin{equation}
\dot{x}= x \left( 1 - \frac{x}{2} + y \right), \quad \dot{y}= y
\left( - 3 + \frac{x}{2} - y \right), \label{0cemo}
\end{equation}
has the irreducible invariant algebraic curves $x=0$, $y=0$ and
$f(x,y)=0$, where $f(x,y):= (x-2)^2 - 2 x y$. Applying the results
described in \cite{seiden4}, it can be shown that this system has
no other irreducible invariant algebraic curve. The function $
V(x,y) = x^{-1/2} y^{1/2} f(x,y)$ is the only algebraic inverse
integrating factor of system (\ref{0cemo}) (modulus multiplication
by non null constants). Since there is no rational inverse
integrating factor, we deduce, by Theorem \ref{0thcggll2}, that
there is no (generalized) Darboux first integral. An elementary
first integral for this system, which is of the form (\ref{0ef}),
is given by:
\[ H(x,y):= \sqrt{2} \sqrt{x} \sqrt{y} + \ln ( x-2 + \sqrt{2} \sqrt{x}
\sqrt{y}) - \ln ( x-2 - \sqrt{2} \sqrt{x} \sqrt{y}). \]
\par
We remark that both Theorems \ref{0thprellesinger} and
\ref{0thcggll2} give a necessary condition to have an elementary
or (generalized) Darboux, respectively, first integral. The
reciprocals to the statements of Theorems \ref{0thprellesinger}
and \ref{0thcggll2} are not true. A result to clarify the easiest
functional class of the first integral once we know the inverse
integrating factor appears in \cite{ChaGaSo06}, see also
\cite{FeLliMa}, where the following theorem is stated:

\begin{theorem} {\rm \cite{ChaGaSo06}} If the system {\rm (\ref{0eq1})} has
a rational inverse integrating factor, then the system has a
(generalized) Darboux first integral. \label{0thfer}
\end{theorem}

In any case, the following Theorem \ref{0thsinger} ensures that
given an algebraic inverse integrating factor, there is a
Liouvillian first integral. The Liouvillian class of functions
contains the rational, algebraic, Darboux and elementary classes
of functions.
\newline

M.F. Singer shows in \cite{Singer} the characterization of the
existence of a Liouvillian first integral for a system
(\ref{0eq1}) by means of its invariant algebraic curves.
\begin{theorem} \index{Singer Theorem} {\sc \cite{Singer}} \
System {\rm (\ref{0eq1})} has a Liouvillian first integral if, and
only if, there is an inverse integrating factor of the form $ V=
\exp \left\{ \int_{(x_0,y_0)}^{(x,y)} \eta \right\} , $ where
$\eta$ is a rational $1$--form such that $d\eta \equiv 0$.
\label{0thsinger}
\end{theorem}
We recall that when $1$--form $\eta$ is such that $d\eta \equiv
0$, we say that it is {\em closed} and if there exists a function
$\varphi$ such that $\eta \, = \, d \varphi$, we say that $\eta$
is {\em exact}. \newline

Taking into account Theorem \ref{0thsinger}, C. Christopher in
\cite{Christopher2} gives the following result, which makes
precise the form of the inverse integrating factor.
\begin{theorem} {\sc \cite{Christopher2}} \ If the system {\rm
(\ref{0eq1})} has an inverse integrating factor of the form $\exp
\left\{ \int_{(x_0,y_0)}^{(x,y)} \eta \right\} , $ where $\eta$ is
a rational $1$--form such that $d\eta \equiv 0$, then there exists
an inverse integrating factor of system {\rm (\ref{0eq1})} of the
form
\[ V= \exp \{D/E \} \prod C_i^{l_i}, \]
where $D$, $E$ and the $C_i$ are polynomials in $x$ and $y$ and
$l_i \in \mathbb{C}$. \label{0thchris}
\end{theorem}
We notice that $C_i=0$ are invariant algebraic curves and $\exp \{
D/E \}$ is an exponential factor for system (\ref{0eq1}). In fact,
since system (\ref{0eq1}) is a real system, we can assume, without
loss of generality, that $V$ is a real function.
\par
Theorem \ref{0thchris} states that the search for Liouvillian
first integrals can be reduced to the search of invariant
algebraic curves and exponential factors. Therefore, if we
characterize the possible cofactors, we have the invariant
algebraic curves of a system and, hence, its Liouvillian or non
Liouvillian integrability. \par Several works study the relation
between the existence of invariant algebraic curves and the
integrability of the system. The existence of an inverse
integrating factor and the functional class it belongs to is
crucial in the resolution of the integrability problem, as
Theorems \ref{0thcggll2}, \ref{0thprellesinger} and \ref{0thchris}
show. A number high enough of invariant algebraic curves of system
(\ref{0eq1}) implies its integrability in one of the rational,
elementary or Liouvillian class, due to Darboux's theorem and
Jouanoulou's improvement, see also Theorems \ref{Main-CGGLl2} and
\ref{Teo-V-KCZ}. The degree of an invariant algebraic curve is not
necessarily related with the integrability class of the system,
see \cite{ChaGra03, Moulin1} and the references therein. \newline

We conclude this part with a theorem that summarizes some
relations between inverse integrating factors and first integrals
of polynomial vector fields.

\begin{theorem} \label{thResum}
Let $\mathcal{X}$ be a planar polynomial vector field.
\begin{itemize}
\item[{\rm (i)}] If $\mathcal{X}$ has a Liouvillian first
integral, then it has a Darboux inverse integrating factor.
\item[{\rm (ii)}] If $\mathcal{X}$ has a Darboux first integral,
then it has a rational inverse integrating factor. \item[{\rm
(iii)}] If $\mathcal{X}$ has a polynomial first integral then it
has a polynomial inverse integrating factor.
\end{itemize}
\end{theorem}
Statement (i) of Theorem \ref{thResum} was proved in \cite{Singer}
and \cite{Christopher2} and statements (ii) and (iii) in
\cite{CGGLl2}.\newline

Another problem related with the inverse integrating factor and
the integrability problem is an inverse problem: given a function
$V(x,y)$, the question is to find (all the) planar differential
systems with $V(x,y)$ as inverse integrating factor. In the case
of searching for a Darboux inverse integrating factor, a very
exhaustive approach to this problem is given in \cite{CLPW08,
LlibPan04, Chara}. The main result of \cite{CLPW08} establishes,
under two generic conditions, {\bf all} the planar polynomial
differential systems with an inverse integrating factor of the
form $V(x,y)=f_1^{\lambda_1}f_2^{\lambda_2} \ldots
f_s^{\lambda_s}, $ where each $f_i(x,y)=0$ is an invariant
algebraic curve of the system and $\lambda_i \in \mathbb{C}$, for
$i=1,2,\ldots,s$, $s \in \mathbb{N}$. We do not reproduce the main
result of \cite{CLPW08} because a lot of notation would need to be
introduced.
\par
In \cite{ChaGiaGin0} another method to construct systems with a
given inverse integrating factor is described. In fact, in 1997
the function $V$ was named {\it null divergence factor}. This
method is a generalization of the classical Darboux method to
generate integrable systems. One of the main results in this paper
is the following one.
\begin{theorem} {\rm \cite{ChaGiaGin0}} \label{thChaGiaGin0}
Let $\mathcal{X}_i = P_i(x,y) \partial_x + Q_i(x,y) \partial_y$, with $i = 1,2, \ldots ,n$, be $\mathcal{C}^1$
vector fields defined in an open subset $\mathcal{U} \subseteq
\mathbb{R}^2$, which have $\mathcal{C}^2$ inverse integrating
factors $V_i(x, y)$, respectively. Then, the
vector field $\mathcal{X} = P(x,y) \partial_x + Q(x,y) \partial_y$ with
\begin{eqnarray*}
P & = & \displaystyle \lambda_0 \frac{\partial V}{\partial y} \, +
\, \sum_{i=1}^{n} \lambda_i \left( \prod_{j=1,j\neq i}^{n}
V_j\right) P_i, \vspace{0.2cm} \\ Q & = & \displaystyle -\lambda_0
\frac{\partial V}{\partial x} \, + \, \sum_{i=1}^{n} \lambda_i
\left( \prod_{j=1,j\neq i}^{n} V_j\right) Q_i,
\end{eqnarray*}
where $\lambda_i$ are arbitrary real numbers for
$i=0,1,2,\ldots,n$, has the inverse integrating factor $V(x,y)$
given by $ \displaystyle V(x,y) \, = \, \prod_{i=1}^{n} V_i(x,y).$
\end{theorem} Indeed, if two systems have the same inverse integrating factor,
a more general system which has such inverse integrating factor
can be constructed, as it is shown in the following proposition.

\begin{proposition} {\rm \cite{ChaGiaGin0}} \label{propChaGiaGin0}
Let $\mathcal{X}_i = P_i(x,y) \partial_x + Q_i(x,y) \partial_y$ with $i=1,2$, be two $\mathcal{C}^1$ vector
fields defined in an open subset $\mathcal{U} \subseteq
\mathbb{R}^2$, which have the same inverse integrating factor
$V(x,y)$. Then, the vector field $\mathcal{X}_1+ \lambda \mathcal{X}_2$ has also the function $V (x, y)$ as an inverse integrating factor, for arbitrary values of the real parameter $\lambda$.
\end{proposition}

This proposition establishes that the set of vector fields with
the same inverse integrating factor forms a $\mathbb{R}$ vector
space. \newline

A polynomial inverse integrating factor allows the study of the
dynamics of system (\ref{0eq1}), because a first integral can be
computed, but it is not so involving as looking for a polynomial
first integral. Indeed, once the degree of a polynomial inverse
integrating factor is fixed, by an ansatz for instance, the
problem of looking for it is reduced to a system of linear
equations on its coefficients. Many authors have used this idea to
find families of planar polynomial differential systems of the
form (\ref{0eq1}) for which all the dynamics can be determined
through an inverse integrating factor. 

In \cite{ChaGiaGin1}, necessary conditions for a planar polynomial vector field to have
a polynomial inverse integrating factor are obtained, see also
\cite{FeLliMa}. All the quadratic systems with a polynomial
inverse integrating factor are determined in \cite{CollFeLli} and
all the quadratic systems with a polynomial first integral are
given in \cite{ChaGLPR}. 
\newline

In \cite{CaLli1} all polynomial first integrals of the non-homogeneous two--dimensional
Lotka--Volterra system of ordinary differential equations are
determined and the role of polynomial inverse integrating factors
is emphasized. Indeed, new first integrals of this class of
systems having a polynomial inverse integrating factor is
presented. The Liouvillian integrability of Lotka-Volterra systems
has been studied in \cite{Moulin, CaGiLli03}.
\newline

In the work \cite{ChGaGi01}, planar differential systems of the
form (\ref{0eq1}) and defined by the sum of homogeneous vector
fields are studied. In particular systems with degenerate infinity
are taken into account. Let us denote by $P_{\rm d}(x,y)$ and
$Q_{\rm d}(x,y)$ the terms of the highest degree ${\rm d}$ in
system (\ref{0eq1}). We say that system (\ref{0eq1}) is {\em of
degenerate infinity} if $x Q_{\rm d}(x,y) - y P_{\rm d}(x,y)
\equiv 0$.  We remark that when a system (\ref{0eq1}) with
degenerate infinity is embedded into a compact space (either by
the Poincar\'e compactification into an sphere or when it is
embedded in the complex projective plane) the line at infinity is
filled with singular points. 

We recall that a real function $H(x,y)$ is said to be $p$-degree homogeneous if $H(\lambda \, x ,
\lambda \, y) \, = \, \lambda^p \, H(x,y)$ for all $(x,y)$ in the
domain of definition of $H(x,y)$ and for all $\lambda \in
\mathbb{R}$, where $p \in \mathbb{Z}$. 
\newline

One of the main results in \cite{ChGaGi01} is the following one.

\begin{theorem} {\rm \cite{ChGaGi01}} \label{thChGaGi01}
Let us consider the following planar polynomial differential
system
\begin{equation} \dot{x} \, = \, P_n(x,y) \, + \, x \, A_{\rm d -1}(x,y), \quad \dot{y} \, = \, Q_n(x,y) \, + \, y \, A_{\rm d
-1}(x,y), \label{eqthChGaGi01+} \end{equation} where $P_n(x,y)$
and $Q_n(x,y)$ are homogeneous real polynomials of degree $n$,
$A_{\rm d -1}(x,y)$ is a real homogeneous polynomial of degree
${\rm d}-1$ and ${\rm d}>n \geq 1$. Let us also consider the
related homogeneous polynomial differential system:
\begin{equation} \dot{x} \, = \, P_n(x,y), \quad \dot{y} \, = \, Q_n(x,y). \label{eqthChGaGi010} \end{equation}
Then, the following statements hold.
\begin{itemize}
\item[{\rm (a)}] If $H(x,y)$ is a $p$-degree homogeneous first integral of
system {\rm (\ref{eqthChGaGi010})}, then $H(x,y)$ is a particular
solution of system {\rm (\ref{eqthChGaGi01+})}.
\item[{\rm (b)}] The homogeneous function $V_{n+1}(x,y) \, := \, x
Q_n(x,y)-yP_n(x,y)$ is an inverse integrating factor of system
{\rm (\ref{eqthChGaGi010})}.
\item[{\rm (c)}] The homogeneous function $V_{n+1}(x,y) \, := \, x
Q_n(x,y)-yP_n(x,y)$ is a particular solution of system {\rm
(\ref{eqthChGaGi01+})}.
\item[{\rm (d)}] If $H(x,y)$ is a $p$-degree homogeneous first integral of
system {\rm (\ref{eqthChGaGi010})}. Then, the function \[ \left(x
Q_n(x,y)-yP_n(x,y)\right) H(x,y)^{\frac{{\rm d}-n}{p}} \] is a
(generalized) Darboux inverse integrating factor of system {\rm
(\ref{eqthChGaGi01+})}.
\end{itemize}
\end{theorem}

The degree of a polynomial inverse integrating factor in relation
with the degree ${\rm d}$ of the system can be bounded under
certain conditions. The conditions established in the following
result come from the embedding of a planar vector field in
$\mathbb{C}P^2$, see \cite{seiden4} for the complete definition of
this embedding. Consider the polynomial differential system
(\ref{0eq1}) with $P$ and $Q$ coprime polynomials of maximum
degree ${\rm d}$. Extending system (\ref{0eq1}) to a differential
equation in the complex projective plane $\mathbb{C}P^2$, a point
$(X_0:Y_0:0) \in \mathbb{C}P^2$ is termed {\it infinite singular
point} of system (\ref{0eq1}) if $(X_0,Y_0) \in \mathbb{C}^2$ is a
root of the homogeneous polynomial $y P_{{\rm d}}(x,y) - x Q_{{\rm
d}}(x,y)$. Here $P_{{\rm d}}$ and $Q_{{\rm d}}$ denote the highest
homogeneous components of $P$ and $Q$ of degree ${\rm d}$.
Following Seidenberg, a singular point is called {\it simple} if
the eigenvalues $\lambda, \mu \in \mathbb{C}$ associated to its
linear part satisfy $\lambda \neq \mu \neq 0$ and $\lambda / \mu
\not\in \mathbb{Q}^+$, where $\mathbb{Q}^+$ stands for the
positive rational numbers. Given a polynomial $F(x,y)$ of degree
$n$, we denote by $\tilde{F}(X,Y,Z)$ its projectivization, that
is, the homogeneous polynomial $\tilde{F}(X,Y,Z) = Z^n
F(X/Z,Y/Z)$. Taking into account these definitions we can state
the following result of Walcher proved in \cite{Walcher}.

\begin{theorem}{\rm \cite{Walcher}} \label{teopoincare2.2}
Let $V(x,y)$ be a polynomial inverse integrating factor of a
polynomial system $\dot{x} = P(x,y)$, $\dot{y} = Q(x,y)$ of degree
${\rm d}$ with $P$ and $Q$ coprime. Assume that the highest
homogeneous components of $P$ and $Q$ of degree ${\rm d}$ are
coprime. If there is an infinite simple critical point of the
system, then the degree of $V$ is exactly ${\rm d}+1$.
\end{theorem}
We observe that Theorem \ref{teopoincare2.2} is also proved in
\cite{ChGa-Poincare} with the additional assumption that the
singularity at infinity $p$ satisfies $\tilde{V}(p) = 0$ where
$\tilde{V}$ is the projectivization of $V$. We remark that the
approach introduced in \cite{Walcher} uses analytical techniques
such as the Poincar\'e--Dulac normal form and the proof given in
\cite{ChGa-Poincare} is completely algebraic and based on the
extension of differential equations to the complex projective
plane and the results of Seidenberg about the reduction of
singularities. The structure of polynomial inverse integrating
factors is also studied by Walcher in \cite{Walcher2,Walcher3}.
\newline

In relation with rational first integrals and in order to state
the main result of \cite{CGGLl2}, we need to introduce some
preliminary concepts, see also \cite{FeLli07}. Let $H = f /g$ be a
rational first integral of a polynomial system (\ref{0eq1}).
According to Poincar\'e \cite{Poin97} we say that $c \in
\mathbb{C} \cup \{\infty \}$ is a {\em remarkable value} of $H$ if
$f +c g$ is a reducible polynomial in $\mathbb{C}[x,y]$ (here,
$c=\infty$ means that $f +c g$ denotes just $g$). In the work
\cite{CGGLl2} it is proved that there are finitely many remarkable
values for a given rational first integral $H$. \par Let now $H$
be a polynomial first integral of degree $n$ of a polynomial
system (\ref{0eq1}). We say that the degree of $H$ is {\it
minimal} between all the degrees of the polynomial first integrals
of (\ref{0eq1}) if any other polynomial first integral of
(\ref{0eq1}) has degree $\geq n$.
\par Assume $H =f / g$ to be a rational first integral. Hence, we
say that $H$ has {\it degree} $n$ if $n$ is the maximum of the
degrees of $f$ and $g$. Moreover, we say that the degree of $H$ is
minimal between all the degrees of the rational first integrals of
system (\ref{0eq1}) if any other rational first integral of
(\ref{0eq1}) has degree $\geq n$. \par Now suppose that $c \in
\mathbb{C}$ is a remarkable value of a rational first integral
$H=f/g$ and that $\prod_{i=1}^r u_i^{\alpha_i}$ is the
factorization of the polynomial $f + c g$ into irreducible factors
in $\mathbb{C}[x,y]$. If some of the $\alpha_i$ is larger than 1,
then we say that $c$ is a {\it critical remarkable value} of $H$
and that $u_i = 0$ having $\alpha_i > 1$ is a critical remarkable
invariant algebraic curve of (\ref{0eq1}) with exponent
$\alpha_i$. \par Finally, let $f$ be a polynomial. We denote by
$\tilde{f}$ the homogeneous part of $f$ of highest degree and this
notation is also used for a Darboux functions like (\ref{0gdf}).
\par The main result of \cite{CGGLl2} is the following one.

\begin{theorem}{\rm \cite{CGGLl2}} \label{Main-CGGLl2}
Suppose that a complex polynomial vector field $\mathcal{X} =
P(x,y) \partial_x + Q(x,y) \partial_y$ of degree ${\rm d}$ with
$P$ and $Q$ coprime has a Darboux first integral $H$ given by {\rm
(\ref{0gdf})} where the polynomials $f_i$ and $g_i$ are
irreducible and the polynomials $g_i$ and $h_i$ are coprime in
$\mathbb{C}[x, y]$. Then the following statements hold.
\begin{itemize}
\item[{\rm (a)}] The inverse integrating factor $V_{\log H}$
associated to the first integral $\log H$ is a rational function,
and it can be written in the form $V_{\log H} = \prod_{i=1}^m
u_i^{k_i}(x,y)$ with $u_i \in \mathbb{C}[x, y]$ irreducible and
$k_i \in \mathbb{Z}$. Moreover, if $\mathcal{X}$ has no rational
first integrals, then $V_{\log H}$ is the unique rational inverse
integrating factor of $\mathcal{X}$.

\item[{\rm (b)}] Assume that $H$ is a minimal polynomial first integral. Then there exists a
polynomial inverse integrating factor.

\item[{\rm (c)}] Suppose that $H =f / g$ is a minimal rational first integral of $\mathcal{X}$ and that $\mathcal{X}$ has
no polynomial first integrals. It is not restrictive to assume
that $f$ and $g$ are irreducible. Then,
\begin{itemize}
\item[{\rm (c.1)}] the rational function
$$
V_{f/g} = \frac{g^2}{\prod_i u_i^{\alpha_i-1}}
$$
where the product runs over all critical remarkable invariant
algebraic curves $u_i = 0$ having exponent $\alpha_i$ is an
inverse integrating factor; and
\item[{\rm (c.2)}] $\mathcal{X}$ has a polynomial inverse integrating factor if and only if $H$ has at
most two critical remarkable values.
\end{itemize}
\end{itemize}
Additionally, if we assume for the first integral {\rm
(\ref{0gdf})} that $f_i$ (respectively $g_j$) are different for $i
=1, \ldots, r$ (respectively $j = 1,\ldots,\ell$), and that it is
complete (i.e. the unique algebraic invariant curves of system
$\mathcal{X}$ are the $f_i = 0$ and the $g_j = 0$), then the
following two statements hold.
\begin{itemize}
\item[{\rm (d)}] If $\mathcal{X}$ has no rational first integrals,
then the inverse integrating factor $V_{\log H}$ associated to the
first integral $\log H$ is the polynomial \[ V_{\log H} = \prod_
{i=1}^r f_i \prod_ {j=1}^\ell g_j^{n_j+1}. \]

\item[{\rm (e)}] If $\tilde{H}$ is a multi--valued function and $\exp(h_j/g_j)$
are exponential factors of $\mathcal{X}$ for $j = 1, \ldots,
\ell$, then $V_{\log H} = \prod_ {i=1}^r f_i \prod_ {j=1}^\ell
g_j^{n_j+1}$ is a polynomial of degree ${\rm d} + 1$.
\end{itemize}
\end{theorem}

In the particular case that $\mu_i = 0$ for $i=1,\ldots, \ell$ in
the expression of (\ref{0gdf}), statement (d) of Theorem
\ref{Main-CGGLl2} can be thought as a generalization of following
result due to Kooij and Christopher \cite{Koo-Ch1} and
independently to \.Zol\c{a}dek \cite{Zol1}.

\begin{theorem} {\rm \cite{Koo-Ch1,Zol1}} \label{Teo-V-KCZ}
Consider a polynomial vector field $\mathcal{X} = P(x,y)
\partial_x + Q(x,y) \partial_y$ of degree ${\rm d}$ with $P, Q \in
\mathbb{C}[x,y]$ (resp. $P, Q \in \mathbb{R}[x,y]$) having $q$
invariant algebraic curves $f_i = 0$ such that the polynomials
$f_i$ are irreducible in $\mathbb{C}[x,y]$ (resp.
$\mathbb{R}[x,y]$ and satisfy that no more than two curves meet at
any point of the plane $\mathbb{C}^2$ (resp. $\mathbb{R}^2$) and
are not tangent at these points, no two curves have a common
factor in their highest order terms and the sum of the degrees of
the curves is ${\rm d}+ 1$. Then, $\prod_{i=1}^r f_i$ is an
inverse integrating factor of $\mathcal{X}$ and $\prod_{i=1}^r
f_i^{\lambda_i}$ for convenient $\lambda_i \in \mathbb{C}$ (resp.
$\lambda_i \in \mathbb{R}$) is a first integral of $\mathcal{X}$.
\end{theorem}

$ $From statement (d) of Theorem \ref{Main-CGGLl2}, the following
result easily follows.

\begin{corollary} {\rm \cite{CGGLl2}} \label{Cor1-CGGLl2}
Suppose that a real polynomial vector field $\mathcal{X} = P(x,y)
\partial_x + Q(x,y) \partial_y$ of degree ${\rm d}$ with $P$ and $Q$
coprime has a Darboux first integral $H$ given by {\rm
(\ref{0gdf})} where the polynomials $f_i$ and $g_j$ are
irreducible, $f_i \neq g_j$, the polynomials $g_j$ and $h_j$ are
coprime in $\mathbb{R}[x,y]$, $\exp(h_j/g_j^{n_j})$ are
exponential factors of $\mathcal{X}$, the $\lambda_i$ and $\mu_j$
are either real numbers, or if some of them is complex then it
appears its conjugate. If $H$ is complete and $\tilde{H}$ is
multi--valued, then $V_{\log H}= \prod_ {i=1}^r f_i \prod_
{j=1}^\ell g_j^{n_j+1}$ is a polynomial of degree ${\rm d} + 1$.
If the system has foci or limit cycles, these are contained in the
set $\{ V_{\log H} \, = \, 0 \}$.
\end{corollary}

Using Corollary \ref{Cor1-CGGLl2} particularized to quadratic
(${\rm d}=2$) polynomial vector fields, in \cite{CGGLl2} it is
obtained the next result.

\begin{corollary} {\rm \cite{CGGLl2}} \label{Cor2-CGGLl2}
Under the assumptions of Corollary {\rm \ref{Cor1-CGGLl2}} there
are no real quadratic polynomial vector fields with a Darboux
first integral {\rm (\ref{0gdf})} and a limit cycle.
\end{corollary}

Some examples of polynomial systems satisfying the assumptions of
Corollary \ref{Cor1-CGGLl2} are the following ones, see again
\cite{CGGLl2}:
\begin{itemize}
\item $\dot{x}=-y - x f_1(x,y)$, $\dot{y}=x - y f_1(x,y)$ where
$f_1(x,y) = x^2+y^2-1 = 0$ is an invariant circle which becomes an
algebraic limit cycle. The origin is a focus and $f_2(x,y)=
x^2+y^2$ is another invariant algebraic curve. The system
possesses the inverse integrating factor $V=f_1 f_2$.

\item The system $\dot{x}=y - 4 x y$, $\dot{y}= -x + x^2 + 2 x y - y^2$
has the invariant algebraic curves $f_1(x,y) = 1-4 x$, $f_2(x,y) =
\sqrt{2} y + (x+y-1) i$ and $f_3(x,y) = \sqrt{2} y - (x+y-1) i$
with $i^2=-1$. The function $V(x,y) = f_1 f_2 f_3$ is an inverse
integrating factor. Notice that the system has a center at $(0,0)$
and a unstable focus at $(1,0)$
\end{itemize}

The recent work \cite{FeLli07} is also devoted to study the
properties of remarkable values. The polynomial $R(x,y):= \prod_i
u_i^{\alpha_i-1}$ defined as the product of all remarkable curves
powered to their respective exponent minus one, is called the {\em
remarkable factor}. From Theorem \ref{Main-CGGLl2}, if $H$ is a
polynomial first integral, then the remarkable factor $R$ is a
polynomial integrating factor of $\mathcal{X}$. Moreover $R$
divides the product $\prod(H+c_i)$ where $c_i$ are all the
critical remarkable values of $H$. Thus the polynomial $V_R =
\prod(H+c_i)/R$ is an inverse integrating factor of the system.
The following theorem gives some relations between the degree of a
system with a polynomial first integral, the degree of its inverse
integrating factor $V_R(x,y)$ and the number of critical
remarkable values.

\begin{theorem} {\rm \cite{FeLli07}} Let $\mathcal{X}$ be a polynomial
vector field of degree ${\rm d}$ and let $H$ be a minimal
polynomial first integral of $\mathcal{X}$. Consider the
remarkable factor $R$ and the polynomial inverse integrating
factor $V_R$. Let $k$ be the number of critical remarkable values.
Then,
\begin{itemize}
\item[{\rm (a)}] $
k(k+{\rm d})  \leq  \deg V_R  =  k \deg H - \deg R  \leq  k(\deg
H-1)  \leq \deg R( \deg R + {\rm d})$ and
\item[{\rm (b)}] $\deg V_R < \deg H$ if and only if $k=1$. Moreover in this case
$\deg V_R = {\rm d}+1$.
\end{itemize}
\label{thFeLLib07}
\end{theorem}

\section{On the center problem}

One of the classical problems in the qualitative theory of planar
analytic differential systems is to characterize the local phase
portrait near an isolated singular point. By using the blow-up
technique, see \cite{Dumor}, this problem can be solved except
when the singularity is monodromic, that is, it is either a focus
or a center. The problem of distinguishing between a center or a
focus is called the {\it center problem}. Another interesting
problem is to know whether there exists or not a local analytic
first integral defined in a neighborhood of a singular point.
These two problems are equivalent when the singularity has
associated nonvanishing complex conjugated eigenvalues. In this
case, translating the singular point at the origin, after a linear
change of variables and a rescaling of the time variable, the
system can be written into the form:
\begin{equation}\label{V-centernondeg}
\dot{x} = -y + f(x,y) \ , \ \dot{y} = x + g(x,y) \ ,
\end{equation}
where $f(x, y)$ and $g(x, y)$ are analytic functions near the
origin without constant nor linear terms. It is well known since
Poincar\'e and Liapunov that system (\ref{V-centernondeg}) has a
center at the origin if and only if there exists a local analytic
first integral of the form $H(x,y) = x^2 + y^2 + F(x, y)$ defined
in a neighborhood of the origin, where $F$ starts with terms of
order higher than 2. We recall here that the {\em
Poincar\'e--Liapunov constants} are the values $V_{2k}$ defined
from the formal power series $H (x, y) = \sum_{n=2}^\infty H_n(x,
y)$, where $H_2(x, y) =(x^2+y^2)/2$ and $H_n(x,y)$ are homogeneous
polynomials of degree $n$ satisfying $\mathcal{X} H =
\sum_{k=2}^\infty V_{2k}(x^2 + y^2)^k$. The origin is a center of
(\ref{V-centernondeg}) if and only if all the Poincar\'e--Liapunov
constants vanish. When $V_{2j} \, = \, 0$ for $j=2,3, \ldots, k-1$
and $V_{2k} \neq 0$, we say that the origin of system
(\ref{V-centernondeg}) is a focus {\em of order $k$}. \par The
existence of invariant algebraic curves is strongly related with
the origin of system (\ref{V-centernondeg}) being a center, as it
is explained in \cite{Christopher1, Schlomiuk1, Schlomiuk2}.
\newline

The proof of the following result is a particular case of a
theorem that was given by Reeb in \cite{Reeb} (see also Mattei and
Moussu \cite{MaMo1} and Moussu \cite{Mo}). For a proof using
elementary methods see \cite{Coll}.

\begin{theorem} {\rm \cite{Reeb}} \label{Teo-V-reeb}
System {\rm (\ref{V-centernondeg})} has a center at the origin if
and only if there is a nonzero analytic inverse integrating factor
in a neighborhood of the origin.
\end{theorem}

In fact, given a system (\ref{V-centernondeg}), the computational
problems of looking for a first integral $H(x,y) = x^2 + y^2 +
\cdots$ or for an inverse integrating factor $V(x,y) = 1 +
\cdots$, where the dots denote terms of higher order, are of the
same difficulty. Thus, the inverse integrating factor offers an
alternative to the solution of the center problem. In
\cite{flows}, it has been noticed that for many systems of type
(\ref{V-centernondeg}) having a center at the origin there is an
inverse integrating factor $V$ with very simple properties which
can be globally defined in all $\mathbb{R}^2$ and which is usually
a polynomial. By contrary, the first integral is, in general, a
complicated expression that can not be written in terms of
elementary functions. \par In particular, when in system
(\ref{V-centernondeg}) the functions $f$ and $g$ are both
quadratic then there exists a polynomial inverse integrating
factor of degree $3$ or $5$, see \cite{Cha94}, whereas the first
integrals are far more complicated, see \cite{lunkevichsibirskii}.
When the functions $f$ and $g$ in system (\ref{V-centernondeg})
are both cubic homogeneous polynomials and the origin is a center,
there exists a polynomial inverse integrating factor of degree at
most $10$, as it is also shown in \cite{Cha94}. In
\cite{GiaNdi96n}, the authors study cubic systems of the form
(\ref{V-centernondeg}) and give some sufficient conditions for the
origin to be a center. This conditions come from the imposition to
the system to have an inverse integrating factor. \par The work
\cite{ChaSab99} is a survey on isochronous centers, that is,
centers of the form (\ref{V-centernondeg}) such that all the
periodic orbits surrounding the origin have the same period. Many
families of isochronous systems are listed and an explicit
expression of an inverse integrating factor is given in each case.
We include here a couple of results of the ones appearing in
\cite{ChaSab99} which we have chosen for being the most known
examples of isochronous centers in the literature. \par The
quadratic systems with a isochronous center at the origin are
characterized in the following result.
\begin{theorem} \label{thChaSab99a} The origin is an isochronous center of a quadratic system
{\rm (\ref{V-centernondeg})} if, and only if, the system can be
brought by means of an affine change of coordinates and a
rescaling of time, to one of the following four systems. For each
case in the list we include the corresponding inverse integrating
factor $V(x,y)$.
\begin{itemize}
\item[{\rm ($S_1$)}] $\dot{x} \, = \, -y + x^2-y^2$, $\dot{y} \, =
\, x(1+2y)$, with $V(x,y)=(1+2y)^2$.
\item[{\rm ($S_2$)}] $\dot{x} \, = \, -y + x^2$, $\dot{y} \, =
\, x(1+y)$, with $V(x,y)=(1+y)^3$.
\item[{\rm ($S_3$)}] $\dot{x} \, = \, -y -\frac{4}{3} x^2$, $\dot{y} \, =
\, x(1-\frac{16}{3}y)$, with $V(x,y)=(3-16y)(9-24y+32x^2)$.
\item[{\rm ($S_4$)}] $\dot{x} \, = \, -y +\frac{16}{3} x^2-\frac{4}{3} y^2$, $\dot{y} \, =
\, x(1+\frac{8}{3}y)$, with $V(x,y)=(3+8y)(9+96y-256x^2+128y^2)$.
\end{itemize}
\end{theorem}
Let us consider a cubic polynomial system of the form
(\ref{V-centernondeg}) and let us assume that it contains no
quadratic terms, that is, it is the sum of a linear system and a
cubic homogeneous system. We say that such a system is cubic and
with homogeneous nonlinearities. The following results
characterizes which of these systems have an isochronous center at
the origin.
\begin{theorem} \label{thChaSab99b} The origin is an isochronous center of a cubic
system with homogeneous nonlinearities of the form {\rm
(\ref{V-centernondeg})} if, and only if, the system can be brought
by means of an affine change of coordinates and a rescaling of
time, to one of the following four systems. For each case in the
list we include the corresponding inverse integrating factor
$V(x,y)$.
\begin{itemize}
\item[{\rm ($S_1^{*}$)}] $\dot{x} \, = \, -y + x^3-3xy^2$, $\dot{y} \, =
\, x+3x^2y-y^3$, with $V(x,y)=(x^2+y^2)^3$.
\item[{\rm ($S_2^{*}$)}] $\dot{x} \, = \, -y + x^3-xy^2$, $\dot{y} \, =
\, x+x^2y-y^3$, with $V(x,y)=(1+2xy)^2$.
\item[{\rm ($S_3^{*}$)}] $\dot{x} \, = \, -y +3x^2y$, $\dot{y} \, =
\, x-2x^3+9xy^2$, with $V(x,y)=(1-3x^2)^4$.
\item[{\rm ($S_4^{*}$)}] $\dot{x} \, = \, -y -3x^2y$, $\dot{y} \, =
\, x+2x^3-9xy^2$, with $V(x,y)=(1+3x^2)^4$.
\end{itemize}
\end{theorem}

In the works \cite{Cha94, Cha95}, Chavarriga writes system
(\ref{V-centernondeg}) in polar coordinates and studies the
existence of inverse integrating factors polynomial in the radial
variable. In \cite{ChaGi96, ChaGi97, ChaGi98}, the authors look
for possible inverse integrating factors for polynomial vector
fields of the form
\begin{equation}\label{V-integ-center}
\mathcal{X} = -y \partial_x + x \partial_y + \mathcal{X}_s \ ,
\end{equation}
where $\mathcal{X}_s$ is a polynomial homogeneous vector field of
degree $s \geq 2$. In particular, they use the quasi--polar
coordinates $(R, \varphi)$ where $R = r^{s-1}$ and $(r, \varphi)$
are the polar coordinates, that is, $x=r \cos\varphi$, $y =r
\sin\varphi$. Next, it is assumed the existence of an inverse
integrating factor $V(R, \varphi)$ of (\ref{V-integ-center}) which
is polynomial in the variable $R$, that is, of the form $V(R ,
\varphi) = \sum_{i=0}^p V_i(\varphi) R^i$ with $V_0(\varphi)
\equiv 1$ and where $V_i(\varphi)$ are homogeneous trigonometrical
polynomials of degree $i (s-1)$. This assumption is clearly
equivalent to impose an inverse integrating factor of the form
$V(R, \varphi) = \prod_{i=1}^p (1+ x_i(\varphi) R)^{\alpha_i}$,
with $\alpha_i \in \mathbb{R}$. The authors try to solve the
system of equations for the unknown functions $x_i(\varphi)$ in
the cases $p=1,2,3$. The case $p=1$ is totally solved. If $p = 2$,
only is solved the case $\alpha_1 = \alpha_2 = (s + 1)/(s - 1) \pm
1 /2$ with arbitrary $s$. Finally, when $p = 3$ the following two
particular cases are investigated: $\alpha_1 = \alpha_2 = \alpha_3
= 5/3$ and either $s = 2$ or $s =3$. \par In \cite{ChaGiGra1},
some invariants are determined from which a formal first integral
for system (\ref{V-integ-center}) can be computed. Moreover, this
technique is applied to the problem of determining the centers of
polynomial vector fields (\ref{V-integ-center}). Recall that a
complete classification of such centers is known when $s = 2, 3$
but only partial results are known in the cases $s = 4$ and $s =
5$. \newline

Theorem \ref{Teo-V-reeb} is used in \cite{Gine2} to find
conditions to have a center. In this work, Gin\'e proposes a
formal power series $V (x, y) = \sum_{n=0}^\infty \bar{V}_n(x,
y)$, where $\bar{V}_0(x, y) \equiv 1$ and $\bar{V}_n(x,y)$ are
homogeneous polynomials of degree $n$ such that $\mathcal{X} V - V
{\rm div} \mathcal{X} = \sum_{k=2}^\infty v_{2k}(x^2 + y^2)^k$,
where $\mathcal{X}$ is given by (\ref{V-integ-center}) and the
constants $v_{2k}$ are called the {\it inverse integrating factor
constants}. Using the above mentioned quasi--polar coordinates
$(R, \varphi)$ it is shown that, if the Poincar\'e--Liapunov
constants $V_k = 0$ for $k = 1, \ldots, m$ and $V_{m+1} \neq 0$,
then $v_{m+1} = -((m+ 1)(s - 1)+ 2) V_{m+1}$. In this spirit, the
paper \cite{Lin} is concerned with the existence of a formal
integrating factor of planar analytic system having a non
degenerate focus or center at the origin and gives an algorithm to
calculate the Poincar\'e--Liapunov constants of any order.
\newline

Given a real analytic planar vector field $\mathcal{X}_0$ with a
center at $p_0 \in \mathbb{R}^2$, in \cite{GiGiLlib1} the authors
say that this center is {\it limit of a linear type center} if
there exists a 1--parameter family $\mathcal{X}_\epsilon$ of
analytic planar vector fields with $\epsilon \geq 0$, defined in a
neighborhood of $p_\epsilon \in \mathbb{R}^2$ and having a non
degenerate center at $p_\epsilon$ for all $\epsilon > 0$
sufficiently small. The main results of \cite{GiGiLlib1} are
summarized as follows.

\begin{theorem}{\rm \cite{GiGiLlib1}} \label{Teo-limitcenters}
Let $\mathcal{X}_0$ be a real analytic planar vector field with a
center at $p_0 \in \mathbb{R}^2$. Then, the following holds:
\begin{itemize}
\item[{\rm (i)}] If $p_0$ is a nilpotent center, then it is limit of a linear type center.
\item[{\rm (ii)}] If $p_0$ is a Hamiltonian degenerate center, then it is limit of a linear type Hamiltonian center.
\item[{\rm (iii)}] If $p_0$ is a time--reversible degenerate center, then it is limit of a time--reversible linear type center.
\end{itemize}
\end{theorem}

In the work \cite{Gine1}, Gin\'e continues the study of the
analytic centers which are limit of linear type centers. It is
proved that if a degenerate center has an analytic inverse
integrating factor $V(x, y)$ which does not vanish near the
center, then this degenerate center is also the limit of a linear
type center (changing the time variable). The idea is as follows.
Assume $V(x,y)$ is an analytic inverse integrating factor of the
analytic vector field $\mathcal{X}_0 = P(x,y) \partial_x + Q(x,y)
\partial_y$ such that $(0,0)$ is a degenerate center and $V(0,0)
\neq 0$. Hence, the rescaled vector field $\mathcal{X}_0 /V$ is
hamiltonian near the origin with analytic first integral $H(x,y)$.
Therefore, since the perturbed vector field $\mathcal{X}_\epsilon
= \mathcal{X}_0 /V + \epsilon (-y \partial_x + x \partial_y)$ is
Hamiltonian too, the origin becomes a linear type center of
$\mathcal{X}_\epsilon$ for all $\epsilon \neq 0$.

\section{Limit cycles}

Let $V(x,y)$ be an inverse integrating factor in the open set
$\mathcal{U} \subset \mathbb{R}^2$ of a
$\mathcal{C}^1(\mathcal{U})$ planar vector field $\mathcal{X}$.
That is, the vector field $\mathcal{X} / V$ has zero divergence,
where defined. If $\mathcal{W}$ is any simply connected component
of $\mathcal{U} \setminus V^{-1}(0)$, then the condition ${\rm
div}(\mathcal{X} / V) \equiv 0$ implies that $\mathcal{X} / V$ is
\emph{Hamiltonian} on $\mathcal{W}$ with $\mathcal{C}^2$
single--valued hamiltonian function $H:\mathcal{W} \to
\mathbb{R}$. Since Hamiltonian systems are area--preserving, hence
have no limit cycles, and $\mathcal{X}$ and $\mathcal{X} / V$ are
topologically equivalent, it follows immediately that, in the
presence of an inverse integrating factor $V$, any limit cycle of
$\mathcal{X}$ lies either in $V^{-1}(0)$ or in a component of
$\mathcal{U} \setminus V^{-1}(0)$ that is not simply connected.
Using the machinery of de Rham cohomology, Giacomini, Llibre, and
Viano eliminated the latter possibility in \cite{GLV}. Hence, they
prove the following theorem.

\begin{theorem} {\rm \cite{GLV}} \label{Teo-vanulaciclo}
Let $\gamma$ be a limit cycle of a $\mathcal{C}^1$ real planar
vector field $\mathcal{X}$ and let $V$ be any inverse integrating
factor of $\mathcal{X}$ defined in some neighborhood of $\gamma$.
Then, $\gamma \subset V^{-1}(0)$.
\end{theorem}
A different proof of Theorem \ref{Teo-vanulaciclo} can be found in
\cite{BerroneGiacomini2, Ga-Sh}. We would like to recall here that
M.V. Dolov in \cite{Dolov0} studies the existence of a single
valued regular integrating factor in a neighborhood of a limit
cycle and presents some connections between an integrating factor
and a limit cycle. Moreover, in the works \cite{Dolov,Dolov2} of
Dolov and coauthors, published before the proof of Theorem
\ref{Teo-vanulaciclo}, it is shown that vector fields with a
Darboux inverse integrating factor of the form $V=\exp(R)$ with
rational $R$ cannot have limit cycles. \newline

Theorem \ref{Teo-vanulaciclo} has been applied in many papers to
study limit cycles of a system as we will see in forthcoming
sections. As an example where this theorem is applied, we would
like to recall the result of Llibre and Rodr\'{\i}guez in
\cite{LliRod04} where it is shown that every finite configuration
of disjoint simple closed curves of the plane is topologically
realizable as the set of limit cycles of a polynomial vector
field. Moreover, the realization can be made by algebraic limit
cycles, and an explicit polynomial vector field exhibiting any
given finite configuration of limit cycles is provided. The proof
of this realization makes use of the inverse integrating factor
and, in particular, of Theorem \ref{Teo-vanulaciclo}. A
generalization of the result of Llibre and Rodr\'{\i}guez is given
in \cite{Peralta05} for systems in higher dimension, that is, it
is shown that any finite configuration of (smooth) cycles in
$\mathbb{R}^n$ can be realized (up to global diffeomorphism) as
hyperbolic and asymptotically stable limit cycles of a polynomial
vector field.
\newline

$ $From Theorem \ref{Teo-vanulaciclo} we have that the knowledge
of an inverse integrating factor for a planar differential system
(\ref{eq1}) implies the knowledge of the number (and location) of
the limit cycles of the system. Many authors have treated the
problem of the existence of an inverse integrating factor. For a
polynomial system (thus defined in the whole $\mathbb{R}^2$), the
knowledge of a polynomial inverse integrating factor solves the
question of the number and location of limit cycles of the
polynomial system, see Section \ref{sectinteg}. In \cite{GaGiSo1},
the authors study the problem of existence of a polynomial inverse
integrating factor in several cases of quadratic vector fields
$\mathcal{X}$. If such an integrating factor $V(x,y)$ exists, then
from Theorem \ref{Teo-vanulaciclo} the curve $V = 0$ is invariant
for $\mathcal{X}$ and any limit cycle of $\mathcal{X}$ lies in
this curve. Therefore, in \cite{GaGiSo1}, the authors study planar
quadratic polynomial vector fields that can have limit cycles and
study the nonexistence of invariant algebraic curves, polynomial
inverse integrating factors and algebraic limit cycles of
arbitrary degree for these systems. Ye Yian-Qian \cite{Yian-Qian}
classified real quadratic systems that can have limit cycles in
the following three families
$$
\dot{x} = \delta x - y + \ell x^2 + M x y + N y^2 \ , \ \
\dot{y}=x (1+a x+ b y) \ ,
$$
according to: family $(I)$ if $a=b=0$; family $(II)$ if $a \neq 0$
and $b=0$; family $(III)$ if $b \neq 0$. In \cite{GaGiSo1} it is proved that
there are not algebraic limit cycles except for $\ell N \delta
\neq 0$ and $M^2-4\ell N \geq 0$ in family $(I)$ (this result is improved in \cite{ChaGaSo} where it is proved that there is
no algebraic limit cycle for family $(I)$). Moreover, they also prove that the polynomial inverse integrating
factors into families $(I)$, $(II)_{N=0}$, $(III)_{a=0}$ and
$(III)_{N=0}$ generically have at most degree $3$. So, in the
studied cases, the existence of polynomial inverse integrating
factor implies the nonexistence of limit cycles or at
most the existence of a circle as a unique limit cycle. 
\newline

Another interesting example of application of Theorem
\ref{Teo-vanulaciclo} is given in the proof of several extensions
to the Bendixson–-Dulac Criterion to study of the number of limit
cycles of planar differential systems, see \cite{GasGia02,
GasGia09, GasGiaLli}. An open set $\mathcal{U} \subseteq
\mathbb{R}^2$ with smooth boundary is said to be $\ell$--connected
if its fundamental group, $\pi_1(\mathcal{U})$ is $Z*
\ldots^{(\ell)} *Z$, or in other words if $\mathcal{U}$ has $\ell$
gaps. The classical Bendixson--Dulac Criterion is the following
proposition, see \cite{GasGia02} for the statement and a short
proof.
\begin{proposition} {\rm (Bendixson--Dulac Criterion)}
\label{propgasgia02} Let $\mathcal{U}$ be an open
$\ell$--connected subset of $\mathbb{R}^2$ with smooth boundary.
Let $\mathcal{X} \, = \, P(x,y) \, \partial_x \, + \,  Q(x,y) \,
\partial_y$ be a vector field of class $\mathcal{C}^1$ in
$\mathcal{U}$. Let $g:\mathcal{U} \to \mathbb{R}$ be a
$\mathcal{C}^1$ function such that
\[ M\, := \, {\rm div}(g\mathcal{X}) \, = \, P \, \frac{\partial
g}{\partial x} \, + \, Q\, \frac{\partial g}{\partial y} \, + \, g
\left( \frac{\partial P}{\partial x} + \frac{\partial Q}{\partial
y}\right)\] does not change sign in $\mathcal{U}$ and vanishes
only on a null measure Lebesgue set, such that $\{ M= 0\} \cap \{
g=0\}$ does not contain periodic orbits of $\mathcal{X}$. Then the
maximum number of periodic orbits of $\mathcal{X}$ contained in
$\mathcal{U}$ is $\ell$. Furthermore, each one of them is a
hyperbolic limit cycle that does not cut $\{g=0\}$ and its
stability is given by the sign of $gM$ over it. \end{proposition}

\section{The zero set of inverse integrating factors}

In Theorem \ref{Teo-vanulaciclo} it is shown that limit cycles
$\gamma$ of a $\mathcal{C}^1$ real planar vector field
$\mathcal{X}$ belong to the zero set of any inverse integrating
factor of $\mathcal{X}$ defined near $\gamma$, that is, $\gamma
\subset V^{-1}(0)$.

In addition to containing any limit cycle of $\mathcal{X}$ lying
in $\mathcal{U}$, the zero set of $V$ is also often connected to
the separatrices of critical points of $\mathcal{X}$ in
$\mathcal{U}$. To understand why, recall that integral curves of
$\mathcal{X}$ that map to themselves under the action of a Lie
group are \emph{invariant solutions} for the Lie group. Recall
that when $\mathcal{Y} = \xi(x,y) \partial_x + \eta(x,y)
\partial_y$ is the infinitesimal generator of a nontrivial local
Lie group of symmetries of $\mathcal{X}$ then the function $V(x,y)
= \det\{\mathcal{X},  \mathcal{Y}\}$ is an inverse integrating
factor of $\mathcal{X}$, as it has already been stated in Section
\ref{sect3}. It is obvious that every solution of $\mathcal{X}$
which remains invariant under the action of the group with
infinitesimal generator $\mathcal{Y}$ must satisfy $V(x,y) = 0$.
In other words, inverse integrating factors must vanish on
invariant solutions.

Based on these ideas, Bluman and Anco \cite{Bluman-Anco} argue
heuristically that separatrices should also lie in $V^{-1}(0)$. Of
course, any saddle loop in a Hamiltonian system is composed of
separatrices not lying in the zero set of the trivial inverse
integrating factor $V \equiv 1$. Nevertheless, the idea has merit
and we expect the zero set of $V$ to play a role in the dynamics
of $\mathcal{X}$ and it is very surprising that this fact was not
completely accomplished until recent times. We repeat verbatim the
following historical development on this issue given in
\cite{BerroneGiacomini}.

\begin{quotation}
``To our best knowledge, J. M. Page was the first author in making
an observation of this kind. Concretely, the idea is developed in
\cite{Pa} of using Lie groups in the computation of singular
solutions to the implicit first order differential equation $ F(x,
y, y') = 0$. The same idea is gathered in pgs. 113 and ss. of
\cite{Pa}, where several examples of calculation of envelopes are
given, and later quoted without variations in \cite{Co}, pgs. 66
and ss.. In pg. 111 of his classical textbook \cite{In}, first
published in 1926, E. L. Ince rescue Page's observation on
envelopes but no other material is added. It would took several
decades until some advance along this line of thinking might be
registered. In this regard, the works of W. H. Steeb, C. E.
Wulfman and G. D. Bluman and S. Kumei must be cited. In
\cite{Steeb}, Steeb discussed the connection between limit cycles
of two-dimensional systems and one-parameter groups of
transformations. In \cite{Wulfman}, Wulfman stated apparently
general conditions on the infinitesimal generator of a Lie group
admitted by a system of autonomous differential equations in order
that an invariant solution is a limit cycle of the system.
However, the argument he offers to support these conditions rests
on a heuristic more than rigorous basis. In turn, chapter 3 of
\cite{Bluman-Kumei} contains a section (Section 3.6) devoted to
discuss the relationships existing between invariant solutions on
one hand and "exceptional paths" on the other. Even though the
developments in this section of the book seems to remain also on a
semi-heuristic level, several examples and exercises are provided
showing how the technique works in particular systems.''
\end{quotation}

In addition to this exhaustive historical description, we also
would like to add the work of Gonz\'alez--Gasc\'on \cite{Gonz-Ga}
where it is pointed out that if there is an infinitesimal
generator of a Lie symmetry $\mathcal{Y}$ of a vector field
$\mathcal{X}$ in $\mathbb{R}^n$, then on the limit cycles
(periodic isolated orbits) of $\mathcal{X}$ it follows that
$\mathcal{X}$ and $\mathcal{Y}$ are parallel. This implies, in the
particular case of planar fields that the associated inverse
integrating factor $V = \det\{ \mathcal{X}, \mathcal{Y} \}$
vanishes on the limit cycle.
\newline

In \cite{BerroneGiacomini} Berrone and Giacomini showed that,
under mild additional hypotheses, the separatrices of
\emph{hyperbolic} saddle--points lying in $\mathcal{U}$ are
contained in $V^{-1}(0)$, and extended this result by showing that
if $\Gamma$ is a compact limit set all of whose critical points
are hyperbolic saddle--points, then under mild conditions $\Gamma
\subset V^{-1}(0)$ holds. Now, we summarize the results in
\cite{BerroneGiacomini}.

It is easy to see that isolated vanishing points of an inverse
integrating factor are singular points of the vector field.
Moreover, for non--degenerate singularities (singularities $p_0
\in \mathbb{R}^2$ of $\mathcal{X}$ with non vanishing Jacobian
determinant $\det (D \mathcal{X} (p_0)) \neq 0$) one has the
following result.

\begin{theorem}{\rm \cite{BerroneGiacomini}} \label{Teo1-V-BG}
Let $p_0$ be a non--degenerate critical point of a $\mathcal{C}^1$
vector field $\mathcal{X}$ and let $V$ be an inverse integrating
factor defined in a neighborhood of $p_0$ and satisfying $V(p_0)
\neq 0$. If $\det (D \mathcal{X} (p_0)) > 0$ then $p_0$ is a
center. On the contrary, when $\det (D \mathcal{X} (p_0)) < 0$,
$p_0$ is a saddle--point.
\end{theorem}

Next theorem is concerned with the stability of isolated zeroes of
an inverse integrating factor.

\begin{theorem}{\rm \cite{BerroneGiacomini}} \label{Teo1-V-BG-2}
Let $p_0$ be an isolated zero of a non--negative inverse
integrating factor $V$ of a $\mathcal{C}^1$ vector field
$\mathcal{X}$ defined in a neighborhood $\mathcal{U}$ of $p_0$.
Then $p_0$ is a stable (resp. unstable) singular point of
$\mathcal{X}$ provided that ${\rm div}  \mathcal{X}|_{\mathcal{U}}
\leq 0$ (resp. $\geq$). Furthermore, $p_0$ is asymptotically
stable (resp. unstable) provided that ${\rm div}
\mathcal{X}|_{\mathcal{U}} < 0$ (resp. $> 0$).
\end{theorem}

When a singularity $p_0$ of $\mathcal{X}$ is a non--isolated zero
of an inverse integrating factor the following result holds. Here,
given an orbit $\gamma_0$ of $\mathcal{X}$, we denote by
$\omega(\gamma_0)$ and $\alpha(\gamma_0)$ its $\omega$--limit set
and $\alpha$--limit set respectively.

\begin{theorem}{\rm \cite{BerroneGiacomini}}\label{Teo1-V-BG-3}
Let $p_0$ be a non--isolated zero of an inverse integrating factor
$V$ of a $\mathcal{C}^1$ vector field $\mathcal{X}$. Then, one of
the following two possibilities may occur:
\begin{itemize}
\item[{\rm (a)}] There exists at least an orbit $\gamma_0$ of
$\mathcal{X}$ (different of $p_0$) such that
$\omega(\gamma_0)=p_0$ or $\alpha(\gamma_0) = p_0$ and
$V|_{\gamma_0} \equiv 0$.

\item[{\rm (b)}] There exists a infinite sequence $\{ \gamma_n
\}_{n \in \mathbb{N}}$ of periodic orbits of $\mathcal{X}$
accumulating at $p_0$ such that $V|_{\gamma_n} \equiv 0$.
\end{itemize}
\end{theorem}

A singularity $p_0$ of the vector field $\mathcal{X}$ is called
{\it strong} if ${\rm div}  \mathcal{X}(p_0) \neq 0$. Otherwise,
when ${\rm div}  \mathcal{X}(p_0) = 0$, it is called {\it weak}.
For a linear strong saddle points, it is easy to see that every
inverse integrating factor must vanish on all four separatrix
curves of the saddle. As it is established by the next theorem,
the situation with nonlinear hyperbolic saddle points of
$\mathcal{C}^1$ systems is entirely analogous to the linear case.
The proof is based on the normal form of $\mathcal{X}$ near the
hyperbolic saddle $p_0$ and the Stable Manifold Theorem.

\begin{theorem}{\rm \cite{BerroneGiacomini}} \label{Teo1-V-BG-4}
Let $p_0$ be a hyperbolic saddle-point of a $\mathcal{C}^1$ vector
field $\mathcal{X}$ and $V$ an inverse integrating factor defined
in a neighborhood $\mathcal{U}$ of $p_0$. Then $V$ vanishes on all
four separatrix curves of the saddle provided that one of the
following conditions holds: (i) $p_0$ is strong; (ii) $p_0$ is
weak and $V(p_0) = 0$.
\end{theorem}

In the work \cite{GaGiGr}, the previous theorem is slightly
improved. If $p_0$ is a hyperbolic saddle point of a
$\mathcal{C}^{k+1}$ vector field $\mathcal{X}$ whose $k^{th}$
saddle quantity is not zero and $V$ is an inverse integrating
factor defined in a neighborhood of $p_0$, then $V(p_0)=0$ (and,
thus, $V$ vanishes on all four separatrix curves of the saddle).
For a full definition of saddle quantities see Subsection
\ref{sect82} in relation with system (\ref{normal-3.0}). \newline

As a corollary of Theorem \ref{Teo1-V-BG-4}, one can ensure the
vanishing of an inverse integrating factor defined near certain
saddle connections. Recall that a saddle connection is a union of
saddle points and orbits connecting them.

\begin{corollary}{\rm \cite{BerroneGiacomini}} \label{Teo1-V-BG-5}
Let $V$ be an inverse integrating factor defined in a region
containing a saddle connection $\Gamma$ whose critical points are
$p_i$ for $i=1,\ldots, n$. If $V$ vanishes at a certain singular
point $p_k$, then $V|_{\Gamma} \equiv 0$.
\end{corollary}

A {\it graphic} $\bar{\Gamma} = \cup_{i=1}^k \phi_i(t) \cup \{
p_1,\ldots,p_k \}$ is formed by $k$ singular points $p_1, \ldots,
p_k$, $p_{k+1} =p_1$ and $k$ oriented regular orbits $\phi_1(t),
\ldots, \phi_k(t)$, connecting them such that $\phi_i(t)$ is an
unstable characteristic orbit of $p_i$ and a stable characteristic
orbit of $p_{i+1}$. A graphic may or may not have associated a
Poincar\'e return map. In case it has one, it is called a {\it
polycycle}.

Now, let us suppose that $\Gamma$ is a {\it graphic}, that is, $\Gamma$ is a limit set which differs from a critical point or a periodic orbit.

\begin{theorem}{\rm \cite{BerroneGiacomini}}\label{Teo1-V-BG-6}
Let $V$ be an inverse integrating factor defined in a region
containing a compact graphic $\Gamma$. Then, the following holds:
\begin{itemize}
\item[{\rm (a)}] $V$ vanishes at a critical point at least of $\Gamma$.

\item[{\rm (b)}] If all the critical points on $\Gamma$ are non--degenerate, then $V|_{\Gamma} \equiv 0$.
\end{itemize}
\end{theorem}

The main results of the paper \cite{Ga-Sh} are generalizations and
extensions of the previous results stated in
\cite{BerroneGiacomini}. A key ingredient in the proof of the
results of \cite{Ga-Sh} is the concept of an integral invariant,
introduced by Poincar\'e in \cite{Poincare} for arbitrary
dimension, and its relation to inverse integrating factors. We
denote by $\phi(t; (x_0,y_0))$ the solution of (\ref{eq1}) passing
through the point $(x_0,y_0) \in \mathcal{U}$ at $t=0$; $\phi(t;
D)$ will denote the image of a domain $D \subset \mathcal{U}$
under the time--$t$ map of the flow generated by the solutions of
system (\ref{eq1}).

\begin{definition} \label{difig1}
Let $\mu : \mathcal{U} \subset \mathbb{R}^2 \to \mathbb{R}$ be a
non--zero integrable function on $\mathcal{U}$. The integral
\begin{equation} \label{ifig2}
\int_{\phi(t; D)} \mu(x,y) \ dx dy
\end{equation}
is an \emph{integral invariant} of system {\rm (\ref{eq1})} if for
any measurable set $D \subset \mathcal{U}$ the integral is
independent of $t$.
\end{definition}

The function $\mu$ is called the \emph{density} of the integral invariant, based
on the obvious hydrodynamic interpretation. Various versions of the following result can be found in textbooks, see for instance \cite{Andronov}. We state it in a form suited to our needs. In \cite{Ga-Sh} is also provided a short proof.

\begin{lemma}{\rm \cite{Poincare}} \label{lifig1}
Let $\mathcal{U}$ be an open subset of $\mathbb{R}^2$, let $V :
\mathcal{U} \to \mathbb{R}$ be a $\mathcal{C}^1$ function, and
define a $\mathcal{C}^1$ function $\mu: \mathcal{U} \setminus
V^{-1}(0) \to \mathbb{R}$ by $\mu = 1 / V$. Then $V$ is an inverse
integrating factor of system {\rm (\ref{eq1})} in $\mathcal{U}$ if
and only if the integral {\rm (\ref{ifig2})} is an integral
invariant for system {\rm (\ref{eq1})} on $\mathcal{U} \setminus
V^{-1}(0)$.
\end{lemma}

By using the relationship between inverse integrating factors and
integral invariants given in Lemma \ref{lifig1}, it is easy to see
the next result. The definition of parabolic or elliptic sector
can be found, for instance, in \cite{DuLlAr}.

\begin{theorem}{\rm \cite{Ga-Sh}} \label{limitpt}
Let $p_0$ be any critical point of system {\rm (\ref{eq1})} at
which there is an elliptic or parabolic sector. If $V$ is any
inverse integrating factor of {\rm (\ref{eq1})} defined on a
neighborhood of $p_0$, then $V(p_0)=0$.
\end{theorem}

In order to state the next result, we recall that a function $f$ is called a {\it Morse
function} if all its critical points are nondegenerate, i.e., the
associated Hessian matrix has maximal rank at all the critical
points. For Morse functions it is well known, see \cite{Hirsch} for
instance, that the set of critical points is discrete, that is,
has no accumulation points.

\begin{theorem}{\rm \cite{Ga-Sh}} \label{limitset}
Let $\Gamma$ be any compact $\alpha$-- or $\omega$--limit set of
system {\rm (\ref{eq1})} that contains a regular point, and let
$V$ be any inverse integrating factor of {\rm (\ref{eq1})} defined
in some neighborhood of $\Gamma$. Depending on the smoothness of
$V$, the following statements hold.
\begin{itemize}
\item[{\rm (a)}] There exists a point $p$ in $\Gamma$ such that $V(p)=0$.
\item[{\rm (b)}] If $V$ is $\mathcal{C}^2$, then either $\Gamma$ contains a point that is an accumulation
      point of isolated critical points of $V$ or $\Gamma \subset V^{-1}(0)$.
\item[{\rm (c)}] If $V$ is real analytic or Morse, then $\Gamma \subset
V^{-1}(0)$.
\end{itemize}
\end{theorem}

Theorem \ref{Teo1-V-BG-4} does not hold, in general, for non--hyperbolic singularities. But it does generalize for
saddle or saddle--node singularities with exactly one non--zero eigenvalue, as the next result shows.

\begin{theorem}{\rm \cite{Ga-Sh}} \label{teoonenuleigen}
Suppose $p_0$ is an isolated singularity of a $\mathcal{C}^1$
vector field $\mathcal{X}$, and that $V$ is an inverse integrating
factor for $\mathcal{X}$ defined in a neighborhood of $p_0$. If
the linear part $D \mathcal{X}(p_0)$ has exactly one zero
eigenvalue, then $V$ vanishes along any separatrix of
$\mathcal{X}$ at $p_0$.
\end{theorem}

We finish this section by stating a corollary of Theorems \ref{Teo1-V-BG-4} and
\ref{teoonenuleigen}.

\begin{corollary}{\rm \cite{Ga-Sh}}
Let $\Gamma$ be a polycycle (or graphic which need not be a limit
set) of system {\rm (\ref{eq1})} and let $V$ be any inverse
integrating factor of {\rm (\ref{eq1})} defined in some
neighborhood of $\Gamma$. Assume that the critical points of {\rm
(\ref{eq1})} that belong to $\Gamma$ are hyperbolic saddles $p_1,
p_2, \ldots, p_n$ or saddles and saddle--nodes $q_1, q_2, \ldots,
q_m$ with exactly one zero eigenvalue. If the separatrices of
$\Gamma$ are such that they always connect either $p_k$ with $p_j$
and $V(p_k)=0$ or $p_k$ with $q_j$ or $q_k$ with $q_j$ then
$\Gamma \subset V^{-1}(0)$.
\end{corollary}

As the authors of \cite{Ga-Sh} remark, the hypothesis in Theorem
\ref{limitset} that $V$ be real analytic does not seem to be
essential. Thus, in \cite{Ga-Sh} it is conjectured that only in
the class $\mathcal{C}^1$ for $V$, Theorem \ref{limitset} remains
valid. This conjecture was solved positively in \cite{EncPer}. In
summary, in that paper the authors prove that there always exists
a smooth inverse integrating factor in a neighborhood of a limit
cycle and obtain a necessary and sufficient condition for the
existence of an analytic one. This condition is expressed in terms
of the Ecalle--Voronin modulus of the associated Poincar\'e map.
We recall that a germ of a map in the set of real analytic
diffeomorphisms near the origin of $\mathbb{R}$ is analytically
embeddable, i.e., it is the time-one map of an analytic vector
field on the line, if and only if its Ecalle--Voronin modulus is
trivial. The embedding properties of the Poincar\'e map are
crucial for the proof of the next theorem.

\begin{theorem}{\rm \cite{EncPer}} \label{Teo-daniel-1}
Let $\gamma$ be a limit cycle of the analytic planar vector field
$\mathcal{X}$. Then there exists a neighborhood $\mathcal{U}$ of
$\gamma$ and a function $V \in \mathcal{C}^\infty(\mathcal{U})$
which is an inverse integrating factor of $\mathcal{X}$ and
vanishes exactly on $\gamma$. Moreover, $V$ can be chosen analytic
if and only if the Ecalle--Voronin modulus of the germ of the
Poincar\'e map of $\mathcal{X}$ along the limit cycle $\gamma$ is
trivial.
\end{theorem}

\begin{corollary}{\rm \cite{EncPer}} \label{Co-daniel-1}
If $\gamma$ is a hyperbolic limit cycle of an analytic vector field $\mathcal{X}$, then
$\mathcal{X}$ admits an analytic inverse integrating factor in a neighborhood of $\gamma$.
\end{corollary}

In addition, in \cite{EncPer} it is also proved that a
$\mathcal{C}^1$ inverse integrating factor of a $\mathcal{C}^1$
planar vector field must vanish identically on the polycycles
which are limit sets of its flow. We recall that a polycycle of a
$\mathcal{C}^1$ vector field is a compact invariant set which
contains both regular and singular points.

\begin{theorem}{\rm \cite{EncPer}} \label{Teo-daniel-2}
Let $\mathcal{X}$ be a $\mathcal{C}^1$ vector field defined in a
domain $\mathcal{U} \subseteq \mathbb{R}^2$. Suppose that $\Gamma
\subset \mathcal{U}$ is a polycycle which is a limit set of
$\mathcal{X}$ and $\mathcal{X}$ has a finite number of singular
points in $\Gamma$. Then if $\mathcal{X}$ admits a $\mathcal{C}^1$
inverse integrating factor $V$ in $\mathcal{U}$, then $\Gamma
\subset V^{-1}(0)$.
\end{theorem}

The main idea of the proof of Theorem \ref{Teo-daniel-2} is to
pull back the vector field $\mathcal{X}$ and the inverse
integrating factor $V$ to the universal cover of $\mathcal{U}
\backslash \{ V \mathcal{X} = 0 \}$ and exploit the fact that
$\mathcal{X}/V$ lifts to a Hamiltonian vector field in the
covering space.

The existence of inverse integrating factors in a neighborhood of
an elementary singularity is also established in \cite{EncPer}.
The regularity of the inverse integrating factor depends on the
kind of singularity and the proof makes crucial use of the theory
of normal forms for planar vector fields. This considerably
extends previous results in \cite{flows}, where the authors prove
for analytic vector fields the existence of a unique analytic
inverse integrating factor in a neighborhood of a strong focus, or
a non--resonant hyperbolic node, or a Siegel hyperbolic saddle.
\par
The following result, which is stated in \cite{Hopf}, is a summary
and a generalization of several results on the existence of a
smooth and non--flat inverse integrating factor $V_0(x,y)$ in a
neighborhood of an isolated singular point, see \cite{flows,
EncPer,GiaVia}.
\begin{theorem} \label{thpoint}
Let the origin be an isolated singular point of {\rm (\ref{eq1})}
and let $\lambda, \mu \in \mathbb{C}$ be the eigenvalues
associated to the linear part of {\rm (\ref{eq1})}. If $\lambda
\neq 0$, then there exists a smooth and non--flat inverse
integrating factor $V(x,y)$ in a neighborhood of the origin.
\end{theorem}

In \cite{EncPer} the existence of an analytic inverse integrating
factor in a neighborhood of a non-degenerate monodromic singular
point of an analytic system is characterized. If the origin is a
non--degenerate center or a strong focus, there exists an analytic
inverse integrating factor. If the origin is a weak focus, by
Theorem \ref{thpoint} we have the existence of a smooth and
non-flat inverse integrating factor, and there exists an analytic
inverse integrating factor if and only if the Ecalle--Voronin
modulus of the associated Poincar\'e map is trivial, see also
Theorem \ref{Teo-daniel-1}. In \cite{EncPer} the first known
examples of real planar analytic vector fields not admitting an
analytic inverse integrating factor in any neighborhood of either
a limit cycle or an isolated singularity are given. \par

In \cite{Hopf}, we show the existence of an inverse integrating
factor in a neighborhood of some degenerate singular points.
\begin{theorem} {\rm \cite{Hopf}} \label{thHopfpunts}
There exists an inverse integrating factor $V(x,y)$, of class at
least $\mathcal{C}^1$, in a neighborhood of the following two
types of singular points: a degenerate focus without
characteristic directions and a nilpotent focus.
\end{theorem}

\section{Bifurcations}

The inverse integrating factor has been shown to be very useful in
many bifurcation problems. The books \cite{HaleKocak, Roussarie}
contain the main concepts and ideas of this theory in the
framework of ordinary differential equations. \par

Consider system (\ref{eq1}) and take a parametric family of
systems of the form
\begin{equation} \label{qq1}
\dot{x} \, = \, \mathcal{P}(x, y, \varepsilon), \quad \dot{y} \, =
\, \mathcal{Q}(x, y, \varepsilon), \end{equation} where
$\mathcal{P}(x, y, \varepsilon)$ and $\mathcal{Q}(x, y,
\varepsilon)$ are analytic functions in $(x,y)$ in the same open
set as $P(x,y)$ and $Q(x,y)$ (or an open set we are interested
in), are analytic for $\varepsilon$ near the origin and coincide
with $P(x, y) $ and $Q(x, y)$ when $\varepsilon = 0$, that is,
$\mathcal{P}(x, y, 0) \, = \, P(x, y) $ and $\mathcal{Q}(x, y, 0)
\, = \, Q (x, y)$. The parameter $\varepsilon$ is called {\em
bifurcation parameter} and we assume it is defined in a
neighborhood of the origin of $\mathbb{R}^k$, with $k \in
\mathbb{N}$; in many cases we consider that $\varepsilon$ is a
real one-dimensional parameter $(k=1)$. For small values of the
norm of $\varepsilon$, we say that the family of systems
(\ref{qq1}) is a perturbation of system (\ref{eq1}). When
$\varepsilon$ takes values near the origin $0<|\varepsilon|<<1$,
the qualitative behavior of system (\ref{qq1}) can change with
respect to the one of system (\ref{eq1}) for $\varepsilon = 0$. In
this case, we say that a bifurcation has occurred. Bifurcation
theory aims at characterizing under which conditions on system
(\ref{eq1}) and its perturbations, this bifurcations eventually
happen and which are their properties. For example, consider a
singular point $p$ of system (\ref{eq1}) and denote by $\lambda$
and $\mu$ the eigenvalues of the linearization of the system
around $p$. If $\lambda, \mu \in \mathbb{R}$ and $\lambda \, \cdot
\, \mu \, <\, 0$, then we say that $p$ is a hyperbolic saddle and
a classical result states that any perturbation of the system in a
neighborhood of this point has the same qualitative behavior, that
is, we have a saddle singular point that, when $\varepsilon$ tends
to zero, tends to $p$. We give the adjective {\em hyperbolic} to
those objects which maintain their qualitative nature under
perturbations. In contrast, if $p$ is a singular point of center
type, i.e. it has a neighborhood filled with periodic orbits, then
a perturbation of the system usually breaks these orbits and the
point can be transformed, for instance, into a singular point of
focus type, i.e. surrounded by orbits that spiral towards (or
from) it. In this case, we say that $p$ is a bifurcation point.
When the considered family (\ref{qq1}) shows all the possible
sample of qualitative behaviors that might occur when perturbing
an object of system (\ref{eq1}), we say that it is an {\em
unfolding}. The minimum number of parameters needed to have an
unfolding is called the {\em codimension}. \par Bifurcation theory
is one of the most current tools used when trying to solve
$16^{th}$ Hilbert problem, part b. This problem was proposed in
1900 by D. Hilbert and asks for the maximum number and possible
configurations of limit cycles that a polynomial system of the
form (\ref{0eq1}) of degree ${\rm d}$ may have, only depending on
the degree ${\rm d}$. For a fixed system, \'Ecalle (1992) and
Il'yashenko (1991) have demonstrated, in a different and
independent way, that the number of limit cycles that the system
may have is finite. However, the problem of determining whether
there exists an upper bound on the number of limit cycles that a
polynomial system of the form (\ref{0eq1}) can have, only
depending on the degree ${\rm d}$ of the system, is still open.
\par As R. Roussarie defines in \cite{Roussarie}, given a
family of systems of the form (\ref{qq1}), a {\em limit periodic
set} is a compact and nonempty subset $\Gamma$ of points so that
there exists some succession $(\varepsilon_n)_n$ which tends to
$\varepsilon_*$ when $n \to + \infty$ such that for every
$\varepsilon_n$, the corresponding system (\ref{qq1}) has a limit
cycle $\gamma_{\varepsilon_n}$ which tends to $\Gamma$, in the
sense of the Haussdorf distance, when $n \to + \infty$. In this
context, it is assumed that the parameters take values in a
compact set. Following an analogous argument to the one used to
prove Poincar\'e--Bendixson Theorem, the structure of limit
periodic sets can be determined. Given a limit periodic set of the
family (\ref{qq1}), we define its {\em cyclicity} as the maximum
number of limit cycles which can be bifurcated from $\Gamma$ in
this family. In \cite{R3}, see also \cite{Roussarie}, R. Roussarie
showed that the existence of a uniform upper bound in the number
of limit cycles of an analytic family (\ref{qq1}) is equivalent to
that each of its limit periodic sets $\Gamma$ has finite
cyclicity. This equivalence and the fact that all the limit
periodic sets in a family (\ref{qq1}) can be determined shows how
bifurcation theory allows to tackle $16^{th}$ Hilbert problem. In
\cite{DRR}, Dumortier, Roussarie and Rousseau established a list
of 121 cases which are all the possible limit periodic sets that
can appear within the family of quadratic systems and proposed a
program, currently unfinished, to study all these graphics to
demonstrate that there is a uniform upper bound for the number of
limit cycles of polynomial systems of degree $2$. \newline

The knowledge of inverse integrating factors for particular
systems has simplified its study and has allowed the understanding
of several bifurcations. In \cite{ChaGiaGin2}, for instance,  the
following family of cubic systems
\[ \begin{array}{lll}
\dot{x} & = & \displaystyle \lambda x - y + \lambda m_1 x^3 + (m_2
- m_1 + m_1 m_2)x^2y+\lambda m_1 m_2 x y^2 + m_2 y^3, \\ \dot{y} &
= & \displaystyle x + \lambda y - x^3 + \lambda m_1 x^2 y + (m_1
m_2 - m_1 -1) x y^2 + \lambda m_1 m_2 y^3, \end{array} \] where
$\lambda$, $m_1$ and $m_2$ are arbitrary real parameters, is
considered. The fact of knowing an inverse integrating factor \[
V(x,y) \, := \, (x^2+y^2)(1+m_1 x^2+ m_1 m_2 y^2) \] for this
family of systems allows the determination of all the bifurcations
within the family. \newline

Indeed, inverse integrating factors allow the understanding of the
bifurcation of limit cycles from many limit periodic sets, as we
explain in this section. The main result used in this context is
Theorem \ref{Teo-vanulaciclo} as it states that any inverse
integrating factor defined in a neighborhood of a limit cycle
needs to vanish on it. We recall that the zero set of a limit
cycle is formed by orbits of the system. \par We split this
section in three subsections depending on the considered limit
periodic sets.

\subsection{Bifurcation from a period annulus}

In this subsection we consider planar differential systems of the
form (\ref{eq1}) with a singular point of center type. The set of
periodic orbits surrounding this point is called its {\em period
annulus}. A perturbation of the system usually breaks these
periodic orbits but some of them might be maintained as limit
cycles for the perturbed system. We say that this periodic orbits
have bifurcated from the period annulus. \par There are several
methods to determine how many limit cycles bifurcate from the
periodic orbits of a period annulus. These methods are based upon
different tools: the Poincar\'e return map, see for instance
\cite{BlowsPerko}; the Poincar\'e–-Pontrjagin-–Melnikov integrals,
see for instance \cite{GH}; averaging theory, see \cite{SV}; and
the inverse integrating factor, see \cite{GLV2, GLV3,VLG}. This
last method also gives the shape of the bifurcated limit cycles up
to any order of the perturbation parameter. \newline

We are going to explain the method described in \cite{GLV2}. Let
us consider a Hamiltonian planar differential system with a center
at the origin:
\[ \dot{x} \, = \, \frac{\partial H}{\partial y}, \quad \dot{y} \,
= \, - \, \frac{\partial H}{\partial x}, \] where $H(x,y)$ is the
Hamiltonian function and it is analytic in a neighborhood of the
origin. We denote by $\mathcal{P}$ the period annulus of the
center at the origin. Any analytic (nonzero) function of the
Hamiltonian is an inverse integrating factor of the system. In
particular, any (nonzero) constant function is an inverse
integrating factor. \par Leu us consider an analytic perturbation
of the previous Hamiltonian system: \begin{equation} \label{hamp}
\dot{x} \, = \, \mathcal{P}(x,y,\varepsilon), \quad \dot{y} \, =
\, \mathcal{Q}(x,y,\varepsilon), \end{equation} where
\[ \mathcal{P}(x,y,\varepsilon) \, : = \,
\frac{\partial H}{\partial y} \, + \, \sum_{k=1}^{\infty}
\varepsilon^k \, f_k(x,y), \quad \mathcal{Q}(x,y,\varepsilon)  \,
:= \, - \, \frac{\partial H}{\partial x} \, + \,
\sum_{k=1}^{\infty} \varepsilon^k \, g_k(x,y), \] and where
$\varepsilon$ is a small real parameter and $f_k(x,y)$, $g_k(x,y)$
are analytic functions in $\mathcal{P} \cup \{(0,0)\}$. Let us
look for an analytic solution
\[ V(x,y,\varepsilon) \, = \, \sum_{k=0}^{\infty} \varepsilon^k \,
V_k(x,y), \] of the partial differential equation
\[ \mathcal{P} \, \frac{\partial V}{\partial x} \, + \, \mathcal{Q}  \, \frac{\partial V}{\partial
y}\, = \, \left( \frac{\partial \mathcal{P}}{\partial x} \, + \,
\frac{\partial \mathcal{Q}}{\partial y} \right) V. \] This partial
differential equation gives a succession of linear differential
equations for the functions $V_k(x,y)$ which can be solved
recursively. \par The equation of order $0$ in $\varepsilon$
implies that $V_0(x,y)$ needs to be an inverse integrating factor
for the unperturbed system. We have, thus, that $V_0 = V_0(H)$ is
a function of the Hamiltonian $H(x,y)$. \par The equation of order
$1$ in $\varepsilon$ gives a linear differential equation in
$V_1(x,y)$ whose non-homogeneous term contains the function
$V_0(x,y)$. Imposing that the function $V$ needs to be periodic
when evaluated on the unperturbed periodic orbits, it can be shown
that $V_0(h)$ needs to be
\[ V_0(h) \, = \, \lambda \, \int_{\{H=h\}} g_1(x,y) \, dx \, - \, f_1(x,y) \,
dy, \] where $\lambda$ is a nonzero real constant, and $\{H=h\}$
denotes the periodic orbit of $\mathcal{P}$ contained in the
$h$--level set of the Hamiltonian $H(x,y)$. Therefore, $V_0(h)$ is
the first Poincar\'e–-Pontrjagin–-Melnikov integral associated to
system (\ref{hamp}). Once we take this expression of $V_0$, we can
solve the linear differential equation for $V_1$ which is
determined up to the sum of an arbitrary function of the
Hamiltonian $H$. \par By induction on $k$ it can be shown that
$V_{k}$ is determined up to the sum of an arbitrary function of
the Hamiltonian $H$ which we denote by $W_k(h)$. When solving the
linear differential equation of order $k+1$ in $\varepsilon$, and
imposing that the function $V$ needs to be periodic on the orbits
of $\mathcal{P}$, it can be shown that $W_k(h)$ corresponds to the
$k+1$ Poincar\'e–-Pontrjagin–-Melnikov integral associated to
system (\ref{hamp}). \par Indeed, in \cite{GLV2}, the authors show
that, fixed a small value of $|\varepsilon|$, the zero sets of the
functions $\sum_{k=0}^{n} \varepsilon^k \, V_k(x,y)$ give
approximations up to order $\varepsilon^n$ of the limit cycles of
system (\ref{hamp}) which bifurcate from $\mathcal{P}$. When
increasing the value of $n$, better approximations of these limit
cycles are obtained and, thus, their shape is determined. \newline

In \cite{GLV3}, this method is generalized to non-Hamiltonian
centers. The paper \cite{VLG} purports a better understanding of
this method as it studies this problem when the first $\ell-1$
Poincar\'e–-Pontrjagin–-Melnikov functions are identically zero.
The main result in this paper is that, in this case, $V_0(h)$ is
the first non identically zero Poincar\'e–-Pontrjagin–-Melnikov
function. \par Most of these ideas are also used in
\cite{GiaLliVia03} to determine semistable limit cycles that
bifurcate from $\mathcal{P}$. Moreover, the method is applied to
study the limit cycles which bifurcate from a Li\'enard system.
\newline

We remark that this method is not only an alternative to the other
methods as it shows how the inverse integrating factor is linked
to bifurcation problems. This method is computationally as
difficult as any other method but, moreover, it provides the shape
of the bifurcated limit cycles.

\subsection{Bifurcation from monodromic $\omega$-limit sets \label{sect82}}

The work \cite{GaGiGr} is concerned with planar real analytic
systems (\ref{eq1}) with an analytic inverse integrating factor
defined in a neighborhood of a regular orbit $\phi(t)$. First of
all it is shown that the inverse integrating factor defines an
ordinary differential equation for the transition map along the
orbit, see equation (\ref{eqvpi}). Taking two transversal sections
$\Sigma_1$ and $\Sigma_2$ based on $\phi(t)$, it is studied the
transition map of the flow of $\mathcal{X}$ in a neighborhood of
$\phi(t)$. This transition map is studied by means of the {\em
Poincar\'e map} $\Pi: \Sigma_1 \to \Sigma_2$. Given a point in
$\Sigma_1$, we consider the orbit of (\ref{eq1}) with it as
initial point and we follow this orbit until it first intersects
$\Sigma_2$.

Let $(\varphi(s), \psi(s)) \in \mathcal{U}$, with $s \in
\mathcal{I} \subseteq \mathbb{R}$ be a parameterization of the
regular orbit $\phi(t)$ between the base points of $\Sigma_1$ and
$\Sigma_2$. Given a point $(x,y)$ in a sufficiently small
neighborhood of the orbit $(\varphi(s), \psi(s))$, we can always
encounter values of the {\em curvilinear coordinates} $(s, n)$
that realize the following change of variables:
$x(s,n)=\varphi(s)- n \psi'(s)$, $y(s,n) =\psi(s) + n
\varphi'(s)$. We remark that the variable $n$ measures the
distance perpendicular to $\phi(t)$ from the point $(x,y)$ and,
therefore, $n=0$ corresponds to the considered regular orbit
$\phi(t)$. We can assume, without loss of generality, that the
transversal section $\Sigma_1$ corresponds to $\Sigma_1 \, := \,
\left\{ s=0 \right\}$ and $\Sigma_2$ to $\Sigma_2 \,  := \,
\left\{ s=L \right\}$, for a certain real number $L>0$. We perform
the change to curvilinear coordinates $(x,y) \mapsto (s,n)$ in a
neighborhood of the regular orbit $n=0$ with $s \in \mathcal{I} \,
= \, [0, L]$. Then, system (\ref{eq1}) is written as the following
ordinary differential equation:
\begin{equation}
\frac{dn}{ds} \, = \, F(s,n) \ . \label{eq2.0}
\end{equation}
We denote by $\Psi(s;n_0)$ the flow associated to the equation (\ref{eq2.0}) with initial condition $\Psi(0;n_0)=n_0$. In these
coordinates, the Poincar\'e map $\Pi: \Sigma_1 \to \Sigma_2$ between these two transversal sections is given by
$\Pi(n_0)=\Psi(L;n_0)$.

We assume the existence of an analytic inverse integrating factor
$V(x,y)$ in a neighborhood of the considered regular orbit
$\phi(t)$ of the analytic system (\ref{eq1}). In fact, when $\Sigma_1 \neq
\Sigma_2$ and no return is involved, there always exists such an inverse integrating factor. The change to curvilinear coordinates gives us an inverse integrating factor for equation (\ref{eq2.0}), denoted by $\tilde{V}(s,n)$ and which satisfies
\begin{equation} \frac{\partial \tilde{V}}{\partial s} \, + \,
\frac{\partial \tilde{V}}{\partial n} \, F(s,n) \, = \,
\frac{\partial F}{\partial n} \, \tilde{V}(s,n). \label{eq3.0}
\end{equation}

Now, we can state one of the main results of \cite{GaGiGr}.

\begin{theorem}{\rm \cite{GaGiGr}}\label{thvpi}
We consider a regular orbit $\phi(t)$ of the analytic system {\rm
(\ref{eq1})} which has an inverse integrating factor $V(x,y)$ of
class $\mathcal{C}^1$ defined in a neighborhood of it and we
consider the Poincar\'e map associated to the regular orbit
between two transversal sections $\Pi: \Sigma_1 \to \Sigma_2$. We
perform the change to curvilinear coordinates and we consider the
ordinary differential equation {\rm (\ref{eq2.0})} with the
inverse integrating factor $\tilde{V}(s,n)$ which is obtained from
$V(x,y)$. In these coordinates, the transversal sections can be
taken such that $\Sigma_1 \, := \, \left\{ s=0\right\}$ and
$\Sigma_2 \, := \, \left\{ s=L\right\}$, for a certain real value
$L>0$. We parameterize $\Sigma_1$ by the real value of the
coordinate $n$. The following identity holds.
\begin{equation}
\tilde{V} \left(L,\Pi(n)\right) \, = \, \tilde{V} \left(0,n\right)
\Pi'(n). \label{eqvpi}
\end{equation}
\end{theorem}

Theorem \ref{thvpi} is the key point to prove Theorems
\ref{th-mult-limv} and \ref{th-mult-loop0}. \par Further, in
\cite{GaGiGr} the authors consider regular orbits whose Poincar\'e
map is a return map and take profit from the result stated in
Theorem \ref{thvpi} in order to study the Poincar\'e map
associated to a limit cycle or to a homoclinic loop, in terms of
the inverse integrating factor. To do that, the following
definition of {\em vanishing multiplicity} of an analytic inverse
integrating factor $V(x,y)$ of the analytic system (\ref{eq1})
over a regular orbit $\phi(t)$ is needed.

\begin{definition}{\rm \cite{GaGiGr}}
Let $V(x,y)$ be an analytic inverse integrating factor of the
analytic system {\rm (\ref{eq1})} and $\phi(t)$ a regular orbit of
it parameterized by $(\varphi(s), \psi(s)) \in \mathcal{U}$, with
$s \in \mathcal{I} \subseteq \mathbb{R}$. Consider the local
change of coordinates $x(s,n)=\varphi(s) - n \psi'(s)$, $y(s,n)
=\psi(s) + n \varphi'(s)$ defined in a neighborhood of the
considered regular orbit $n=0$ and take the following Taylor
development around $n=0$:
\begin{equation} \label{eq-vnonula.0}
V(x(s,n), y(s,n)) \, = \, n^{m} \, v(s) \, + \, O(n^{m+1}) ,
\end{equation}
where $m$ is an integer with $m \geq 0$ and the function $v(s)$ is
not identically null, we say that $V$ has multiplicity $m$ on $\phi(t)$.
\end{definition}

In fact, in \cite{GaGiGr} it is proved that $v(s) \neq 0$ for any
$s \in \mathcal{I}$, and thus, the vanishing multiplicity of $V$
on $\phi(t)$ is well--defined over all its points. \newline

Let us consider as regular orbit a limit cycle $\gamma$ and  we
use the parameterization of $\gamma$ in curvilinear coordinates
$(s,n)$ with $s \in [0, L)$. Thus, the Poincar\'e map associated
to $\gamma$ is $\Pi(n_0) \, = \, \Psi(L;n_0)$. It is well known
that $\Pi$ is analytic in a neighborhood of $n_0=0$. We recall
that the periodic orbit $\gamma$ is a limit cycle if, and only if,
the Poincar\'e return map $\Pi$ is not the identity. If $\Pi$ is
the identity, we have that $\gamma$ belongs to a period annulus.
We recall the definition of multiplicity of a limit cycle:
$\gamma$ is said to be a limit cycle of {\em multiplicity} $1$ if
$\Pi'(0) \neq 1$ and $\gamma$ is said to be a limit cycle of
multiplicity $m$ with $m \geq 2$ if $\Pi(n_0) \, = \, n_0 \, + \,
\beta_m \, n_0^m \, + \, O(n_0^{m+1})$ with $\beta_m \neq 0$.
Then, one has the following result for limit cycles.

\begin{theorem}{\rm \cite{GaGiGr}} \label{th-mult-limv}
Let $\gamma$ be a periodic orbit of the analytic system {\rm
(\ref{eq1})} and let $V$ be an analytic inverse integrating factor
defined in a neighborhood of $\gamma$.
\begin{itemize}
\item[{\rm (a)}] If $\gamma$ is a limit cycle of multiplicity $m$,
then $V$ has vanishing multiplicity $m$ on $\gamma$.

\item[{\rm (b)}] If $V$ has vanishing multiplicity $m$ on $\gamma$, then $\gamma$ is a limit cycle of multiplicity $m$ or it
belongs to a continuum of periodic orbits.
\end{itemize}
\end{theorem}

Since the Poincar\'e map of a periodic orbit is an analytic function and the multiplicity of a
limit cycle is a natural number, the following corollary is obtained.

\begin{corollary}{\rm \cite{GaGiGr}} \label{cor-mult-limv}
Let $\gamma$ be a periodic orbit of the analytic system {\rm
(\ref{eq1})} and let $V$ be an inverse integrating factor of class
$\mathcal{C}^1$ defined in a neighborhood of $\gamma$. We take the
change to curvilinear coordinates $x(s,n)=\varphi(s) - n
\psi'(s)$, $y(s,n) =\psi(s) + n \varphi'(s)$ defined in a
neighborhood of $\gamma$. If we have that the leading term in the
following development around $n=0$:
\[
V(x(s,n), y(s,n)) \, = \, n^{\rho} \, v(s) \, + \, o(n^{\rho}) ,
\]
where $v(s) \not\equiv 0$ is such that either $\rho=0$ or $\rho
> 1$ and $\rho$ is not a natural number, then $\gamma$ belongs
to a continuum of periodic orbits.
\end{corollary}

A regular orbit $\phi(t) = (x(t), y(t))$ of (\ref{eq1}) is
called a {\it homoclinic orbit} if $\phi(t) \to p_0$ as $t \to \pm
\infty$ for some singular point $p_0$. A {\it homoclinic loop}
is the union $\Gamma = \phi(t) \cup \{p_0\}$. We assume that $p_0$
is a hyperbolic saddle, that is, a critical point of system
(\ref{eq1}) such that the eigenvalues of the Jacobian matrix $D
\mathcal{X}(p_0)$ are both real, different from zero and of
contrary sign. We remark that this type of graphics always has
associated (maybe only its inner or outer neighborhood) a
Poincar\'e return map $\Pi : \Sigma \to \Sigma$ with $\Sigma$ any
local transversal section through a regular point of $\Gamma$. We will assume that $\Gamma$ is a compact invariant set. A goal in \cite{GaGiGr} is to study the cyclicity of the described
homoclinic loop $\Gamma$ in terms of the vanishing multiplicity of
an inverse integrating factor. Roughly speaking, the {\it
cyclicity} of $\Gamma$ is the maximum number of limit cycles which
bifurcate from it under an analytic perturbation of the analytic system (\ref{eq1}). Before state the result for homoclinic loops, we recall briefly that the first saddle quantity is $\alpha_1 \, = \, {\rm div\, } \mathcal{X}(p_0)$ and
it classifies the point $p_0$ between being strong (when $\alpha_1
\neq 0$) or weak (when $\alpha_1=0$). If $p_0$ is a weak saddle
point, the saddle quantities are the obstructions for it to be
analytically orbitally linearizable. In order to define the next saddle quantities associated to $p_0$, we translate the saddle--point $p_0$ to the origin of coordinates and we make a linear change of variables so that its
unstable (resp. stable) separatrix has the horizontal (resp.
vertical) direction at the origin. Let $p_0$ be a weak hyperbolic
saddle point situated at the origin of coordinates and whose
associated eigenvalues are taken to be $\pm 1$ by a rescaling of
time, if necessary. Then, it is well known the existence of an analytic near--identity change
of coordinates that brings the system into:
\begin{equation}\label{normal-3.0}
\begin{array}{lll}
\dot{x} & = & \displaystyle  x \, + \, \sum_{i=1}^{k-1} a_i \,
x^{i+1} y^i \, + \, a_k \, x^{k+1} y^k \, + \, \cdots \ ,
\vspace{0.2cm} \\ \dot{y} & = & \displaystyle -y \, - \,
\sum_{i=1}^{k-1} a_i \, x^{i} y^{i+1} \, - \, b_k \, x^k y^{k+1}
\, + \, \cdots \ , \end{array}
\end{equation}
with $a_k - b_k \neq 0$ and where the dots denote terms of higher
order. The first non--vanishing saddle quantity is defined by
$\alpha_{k+1} := a_k-b_k$.

\begin{theorem}{\rm \cite{GaGiGr}} \label{th-mult-loop0}
Let $\Gamma$ be a compact homoclinic loop through the hyperbolic
saddle $p_0$ of the analytic system {\rm (\ref{eq1})} whose
Poincar\'e return map is not the identity. Let $V$ be an analytic
inverse integrating factor defined in a neighborhood of $\Gamma$
with vanishing multiplicity $m$ over $\Gamma$. Then, $m \geq 1$
and the first possible non--vanishing saddle quantity is
$\alpha_m$. Moreover,
\begin{itemize}
\item[{\rm (i)}] the cyclicity of $\Gamma$ is $2m-1$, if
$\alpha_m \, \neq \, 0$,

\item[{\rm (ii)}] the cyclicity of $\Gamma$ is $2m$, otherwise.
\end{itemize}
\end{theorem}

In addition, in \cite{GaGiGr} it is described one obstruction to
the existence of an analytic inverse integrating factor defined in
a neighborhood of certain homoclinic loops. First of all, we
recall some concepts. By an affine change of coordinates, in a
neighborhood of a hyperbolic saddle, any analytic system can be
written as $ \dot{x} = \lambda x + f(x,y)$, $\dot{y} = \mu y +
g(x,y)$, where $f$ and $g$ are analytic in a neighborhood of the
origin with lowest terms at least of second order and
$\mu<0<\lambda$. This hyperbolic saddle is analytically {\em
orbitally linearizable} if there exists an analytic near--identity
change of coordinates transforming the system to $\dot{x} =
\lambda x h(x,y)$, $\dot{y} =  \mu y h(x,y)$ with $h(0,0)=1$. On
the other hand, when $\mu/\lambda = -q/p \in \mathbb{Q}^{-}$ with
$p$ and $q$ natural and coprime numbers, the saddle is called
$p:q$ resonant.

\begin{proposition}{\rm \cite{GaGiGr}} \label{prop-noV0}
Suppose that the analytic system {\rm (\ref{eq1})} has a
homoclinic loop $\Gamma$ through the hyperbolic saddle point $p_0$
which is not orbitally linearizable, $p:q$ resonant and strong $(p
\neq q)$. Then, there is no analytic inverse integrating factor
$V(x,y)$ defined in a neighborhood of $\Gamma$.
\end{proposition}

\subsection{Generalized Hopf Bifurcation}

Let us consider a planar real system (\ref{eq1}), $\dot{x} =
P(x,y)$, $\dot{y} = Q(x,y)$ and suppose that it is analytic near
an isolated monodromic singular point $p_0$ which we assume to be
at the origin. We associate to system (\ref{eq1}) the vector field
$\mathcal{X}_0 \, = \, P(x,y) \partial_x \, + \, Q(x,y)
\partial_y$. We consider an analytic perturbation of system
(\ref{eq1}) of the form:
\begin{equation}
\dot{x} \, = \, P(x,y) \, + \, \bar{P}(x,y,\varepsilon), \qquad
\dot{y} \, = \, Q(x,y) \, + \, \bar{Q}(x,y,\varepsilon),
\label{eq2}
\end{equation} where $\varepsilon \in \mathbb{R}^p$ is the perturbation parameter,
$0<\|\varepsilon\|<<1$ and the functions
$\bar{P}(x,y,\varepsilon)$ and $\bar{Q}(x,y,\varepsilon)$ are
analytic for $(x,y) \in \mathcal{U}$, analytic in a neighborhood
of $\varepsilon=0$ and $\bar{P}(x,y,0)= \bar{Q}(x,y,0)\equiv 0$.
We associate to this perturbed system (\ref{eq2}) the vector field
$\mathcal{X}_\varepsilon \, = \, (P(x,y) \, + \,
\bar{P}(x,y,\varepsilon))\partial_x \, + \, (Q(x,y) \, + \,
\bar{Q}(x,y,\varepsilon)) \partial_y$.
\par We say that a limit cycle $\gamma_\varepsilon$ of system
(\ref{eq2}) {\em bifurcates from the origin} if it tends to the
origin (in the Hausdorff distance) as $\varepsilon \to 0$. We are
interested in giving a sharp upper bound for the number of limit
cycles which can bifurcate from the origin $p_0$ of system
(\ref{eq1}) under any analytic perturbation with a finite number
$p$ of parameters. The word sharp means that there exists a system
of the form (\ref{eq2}) with exactly that number of limit cycles
bifurcating from the origin, that is, the upper bound is
realizable. This sharp upper bound is called the {\em cyclicity}
of the origin $p_0$ of system (\ref{eq1}) and will be denoted by
${\rm Cycl}(\mathcal{X}_\varepsilon,p_0)$ all along this section.

In \cite{Hopf} we consider systems of the form (\ref{eq1}) where
the origin $p_0$ is a focus singular point of the following three
types: non-degenerate, degenerate without characteristic
directions and nilpotent. The results of \cite{Hopf} do not
establish that the cyclicity of this type of singular points is
finite but give an effective procedure to study it. In the three
mentioned types of focus points, we will consider a change to
(generalized) polar coordinates which embed the neighborhood
$\mathcal{U}$ of the origin into a cylinder $C\, = \, \left\{
(r,\theta) \in \mathbb{R} \times \mathcal{S}^1 \, : \, |r|<\delta
\right\}$ for a certain sufficiently small value of $\delta>0$.
This change to polar coordinates is a diffeomorphism in
$\mathcal{U} - \{(0,0)\}$ and transforms the origin of coordinates
to the circle of equation $r=0$. In these new coordinates, system
(\ref{eq1}) can be seen as a differential equation over the
cylinder $C$ of the form:
\begin{equation}
\frac{dr}{d\theta} \, = \, \mathcal{F}(r,\theta), \label{eq3}
\end{equation} where $\mathcal{F}(r,\theta)$ is an analytic
function in $C$. We consider an inverse integrating factor $V(r,\theta)$ of equation
(\ref{eq3}), that is, a function $V: C \to \mathbb{R}$ of class
$\mathcal{C}^1(C)$, which is non locally null and which satisfies
the following partial differential equation:
\[ \frac{\partial V(r,\theta)}{\partial \theta} \, + \, \frac{\partial
V(r,\theta)}{\partial r} \, \mathcal{F}(r,\theta) \, = \,
\frac{\partial \mathcal{F}(r,\theta)}{\partial r} \, V(r,\theta).
\] We remark that since $V(r,\theta)$ is a function defined over
the cylinder $C$ it needs to be $T$--periodic in $\theta$, where
$T$ is the minimal positive period of the variable $\theta$, that
is, we consider the circle $\mathcal{S}^1 \, = \, \mathbb{R} /
[0,T]$. The function $V(r,\theta)$ is smooth
($\mathcal{C}^{\infty}$) and non--flat in $r$ in a neighborhood of
$r=0$.

Let us consider the Taylor expansion of the function $V(r,\theta)$
around $r=0$: $ V(r,\theta) \, = \, v_m(\theta) \, r^m \, + \,
\mathcal{O}(r^{m+1}), $ where $v_m(\theta) \not\equiv 0$ for
$\theta \in \mathcal{S}^1$ and $m$ is an integer number with $m
\geq 0$. We say that $m$ is the {\em vanishing multiplicity}
of $V(r,\theta)$ on $r=0$. The uniqueness of
$V(r,\theta)$ implies that the number $m$ corresponding to the
vanishing multiplicity of $V(r,\theta)$ on $r=0$ is well--defined.
\newline

We consider a system (\ref{eq1}) of the
form:
\begin{equation} \label{eqdeg}
\dot{x} \, = \, p_d(x,y)\, +\, P_{d+1}(x,y), \qquad \dot{y} \, =
\, q_d(x,y)\, +\, Q_{d+1}(x,y),
\end{equation}
where $d  \geq 1$ is an odd number, $p_d(x,y)$ and $q_d(x,y)$ are
homogeneous polynomials of degree $d$ and $P_{d+1}(x,y),
Q_{d+1}(x,y) \in \mathcal{O}(\|(x,y)\|^{d+1})$. We assume that $p_d^2(x,y)+q_d^2(x,y) \not\equiv 0$. A {\em characteristic direction} for the origin of system (\ref{eqdeg}) is a linear factor in $\mathbb{R}[x,y]$ of the
homogeneous polynomial $xq_d(x,y)-yp_d(x,y)$. If there are no characteristic
directions, then the origin is a monodromic singular point of
system (\ref{eqdeg}), that is, it is either a center or a focus.

In relation with system (\ref{eq2}), an analytic  perturbation field
$(\bar{P}(x,y,\varepsilon), \bar{Q}(x,y,\varepsilon))$ is said to
have subdegree $s$ if $(\bar{P}(x,y,\varepsilon), \bar{Q}(x, y,
\varepsilon)) = \mathcal{O}( \|(x,y)\|^s )$. In this case, we
denote by $\mathcal{X}_\varepsilon^{[s]}\, = \, (P(x,y) \, + \,
\bar{P}(x,y,\varepsilon))\partial_x \, + \, (Q(x,y) \, + \,
\bar{Q}(x,y,\varepsilon)) \partial_y$ the vector field associated
to such a perturbation.

\begin{theorem}{\rm \cite{Hopf}} \label{thdeg}
We assume that the origin $p_0$ of system {\rm (\ref{eqdeg})} is
monodromic and without characteristic directions. Take polar coordinates $x= r \cos\theta$, $y= r \sin\theta$ and let
$V(r,\theta)$ be an inverse integrating factor of the
corresponding equation {\rm (\ref{eq3})} which has a Laurent
expansion in a neighborhood of $r=0$ of the form $ V(r, \theta) \,
= \, v_{m}(\theta) \, r^m \, + \, \mathcal{O}(r^{m+1}), $ with
$v_m(\theta) \not\equiv 0$ and $m \in \mathbb{Z}$.
\begin{itemize}
\item[{\rm (i)}] If $m \leq 0$ or $m$ is even, then the origin of
system {\rm (\ref{eqdeg})} is a center. \item[{\rm (ii)}] If the
origin of system {\rm (\ref{eqdeg})} is a focus, then $m \geq 1$,
$m$ is an odd number and the cyclicity ${\rm
Cycl}(\mathcal{X}_\varepsilon,p_0)$ of the origin of system {\rm
(\ref{eqdeg})} satisfies ${\rm Cycl}(\mathcal{X}_\varepsilon,p_0)
\geq  (m+d)/2-1$. In this case $m$ is the vanishing multiplicity
of $V(r,\theta)$ on $r=0$.
\begin{itemize}
\item[{\rm (ii.1)}] If, moreover, the focus is non--degenerate $(d=1)$,
then the aforementioned lower bound is sharp, that is, ${\rm
Cycl}(\mathcal{X}_\varepsilon,p_0) =  (m-1)/2$.
\item[{\rm (ii.2)}] If only perturbations whose subdegree is greater than or equal to $
d$ are considered, then the maximum number of limit cycles which
bifurcate from the origin is $(m-1)/2$, that is, ${\rm
Cycl}(\mathcal{X}_\varepsilon^{[d]},p_0) =  (m-1)/2$.
\end{itemize}
\end{itemize}
\end{theorem}

\begin{remark} \label{remdeg} From the proof of Theorem {\rm \ref{thdeg}},
it follows that if there exists an inverse integrating factor
$V_0(x,y)$ of system {\rm (\ref{eqdeg})} such that $V_0(r \cos
\theta, \, r  \sin \theta)/r^d$ has a Laurent expansion in a
neighborhood of $r=0$, then the exponents of the leading terms of
$V_0(r \cos \theta, \, r  \sin \theta)/r^d$ and $V(r,\theta)$
coincide. Thus, the vanishing multiplicity $m$ can be computed
without passing the system to polar coordinates.
\end{remark}

We assume that the origin of system (\ref{eqdeg}) is a focus
without characteristic directions and that the vanishing
multiplicity of an inverse integrating factor on it is $m$. If
system (\ref{eqdeg}) is written as $\dot{x} \, = \, P(x,y)$ and
$\dot{y} \, = \, Q(x,y)$, then the system:
\begin{equation} \label{eqdegp1} \dot{x} \, = \, P(x,y) \, + \, x
\, K(x,y,\varepsilon), \qquad \dot{y} \, = \, Q(x,y)  \, + \, y \,
K(x,y,\varepsilon),
\end{equation} where $L \,
: =\, (m+d)/2-1$ and $\displaystyle K(x,y,\varepsilon) \, = \,
\sum_{i=0}^{k-1} \varepsilon^{k-i} \, a_i \, (x^2+y^2)^{i+
\frac{d-1}{2}}$, has at least $(m+d)/2-1$ limit cycles bifurcating
from the origin for convenient values of the real parameters
$a_i$. We recall that both $m$ and $d$ are odd and $d\geq 1$, $m
\geq 1$. \newline

We say that the origin of system (\ref{eq1}) is a {\em nilpotent
singular point} if it is a degenerate singularity that can be
written as:
\begin{equation}\label{eqnil}
\dot{x} = y + P_2(x,y) \ , \ \dot{y} = Q_2(x,y) \ ,
\end{equation}
with $P_2(x,y)$ and $Q_2(x,y)$ analytic functions near the origin
without constant nor linear terms. The following theorem is due to
Andreev \cite{Andreev} and it solves the monodromy problem for the
origin of system (\ref{eqnil}).

\begin{theorem} \label{thandreev} {\sc \cite{Andreev}} Let $y=F(x)$ be the solution of $y + P_2(x,y)
= 0$ passing through $(0, 0)$. Define the functions $f(x) = Q_2(x,
F(x)) = a x^\alpha + \cdots$ with $a \neq 0$ and $\alpha \geq 2$
and $\phi(x) = (\partial P_2/\partial x \, + \, \partial Q_2 /
\partial y)(x,F(x))$. We have that either
$\phi(x) = b x^\beta + \cdots$ with $b \neq 0$ and $\beta \geq 1$
or $\phi(x) \equiv 0$. Then, the origin of {\rm (\ref{eqnil})} is
monodromic if, and only if, $a < 0$, $\alpha = 2 n-1$ is an odd
integer and one of the following conditions holds: (i) $\beta > n-1$; (ii) $\beta = n-1$ and $b^2+4 a n < 0$; (iii) $\phi(x) \equiv 0$.
\end{theorem}

\begin{definition} \label{defnil} We consider a system of the form {\rm (\ref{eqnil})} with the
origin as a monodromic singular point. We define its {\em Andreev
number} $n \geq 2$ as the corresponding integer value given in
Theorem {\em \ref{thandreev}}.
\end{definition}

We consider system (\ref{eqnil}) and we assume that the origin is
a nilpotent monodromic singular point with Andreev number $n$.
Then, the change of variables
\begin{equation}\label{change1}
(x,y) \, \mapsto \, (x, y-F(x)),
\end{equation}
where $F(x)$ is defined in Theorem \ref{thandreev}, and the
scaling
\begin{equation}\label{change2}
(x,y) \, \mapsto \, (\xi \, x, -\xi \, y),
\end{equation}
with $\xi = (-1/a)^{1/(2-2n)}$, brings system (\ref{eqnil}) into
the following analytic form for monodromic nilpotent singularities
\begin{equation}\label{eqnil2}
\dot{x}\, = \, y\, (-1 + X_1(x,y)), \quad \dot{y} \, =\, f(x) + y
\, \phi(x) + y^2\, Y_0(x,y),
\end{equation}
where $X_1(0,0)=0$, $f(x) = x^{2n-1} + \cdots$ with $n \geq 2$ and
either $\phi(x) \equiv 0$ or $\phi(x) = b x^\beta + \cdots$ with
$\beta \geq n-1$. We remark that we have relabelled the functions
$f(x)$, $\phi(x)$ and the constant $b$ with respect to the ones
corresponding to system (\ref{eqnil}).

We are going to transform system (\ref{eqnil2}) to an equation
over a cylinder of the form (\ref{eq3}). The transformation
depends on the Andreev number $n$ and it is given through the {\it
generalized trigonometric functions} defined by Lyapunov as the unique solution $x(\theta) = {\rm Cs} \,
\theta$ and $y(\theta) = {\rm Sn}\, \theta$ of the following
Cauchy problem
\begin{equation}\label{Cauchy}
\frac{d x}{d \theta} \, =\, -y , \ \frac{d y}{d \theta}\, =\,
x^{2n-1}, \qquad x(0)=1 , \, y(0)=0  .
\end{equation}
We introduce in $\mathbb{R}^2 \backslash \{(0,0)\}$ the change to
{\it generalized polar coordinates}, $\,  (x,y) \mapsto
(r,\theta)$, defined by  \begin{equation} \label{change} x\, =\, r
\, {\rm Cs}\, \theta , \qquad y\, =\, r^n \, {\rm Sn}\, \theta.
\end{equation}

We consider the following definition, which will be used in the following Theorem
\ref{thnil}.

\begin{definition} \label{defnilqh}
An analytic perturbation vector field $(\bar{P}(x,y,\varepsilon),
\bar{Q}(x,y,\varepsilon))$ is said to be {\em
$(1,n)$--quasihomogeneous of weighted subdegrees $(w_x,w_y)$} if
$\bar{P}(\lambda x,\lambda^n y,\varepsilon) \, = \,
\mathcal{O}(\lambda^{w_x})$ and $\bar{Q}(\lambda x,\lambda^n
y,\varepsilon) \, = \, \mathcal{O}(\lambda^{w_y})$. In this case,
we denote by $\mathcal{X}_\varepsilon^{[w_x,w_y]}\, = \, (P(x,y)
\, + \, \bar{P}(x,y,\varepsilon))\partial_x \, + \, (Q(x,y) \, +
\, \bar{Q}(x,y,\varepsilon)) \partial_y$ the vector field
associated to such a perturbation.
\end{definition}

The following theorem is one of the main results of \cite{Hopf}.
The symbol $\lfloor x \rfloor$ denotes the integer part of $x$.

\begin{theorem}{\rm \cite{Hopf}} \label{thnil}
We assume that the origin of system {\rm (\ref{eqnil})} is
monodromic with Andreev number $n$. Take generalized polar coordinates (\ref{change}) and let $V(r,\theta)$ be an
inverse integrating factor of the corresponding equation {\rm
(\ref{eq3})} which has a Laurent expansion in a neighborhood of
$r=0$ of the form $ V(r, \theta) \, = \, v_{m}(\theta) \, r^m \, +
\, \mathcal{O}(r^{m+1}), $ with $v_m(\theta) \not\equiv 0$ and $m
\in \mathbb{Z}$.
\begin{itemize}
\item[{\rm (i)}] If $m \leq 0$ or $m+n$ is odd, then the origin of
system {\rm (\ref{eqnil})} is a center. \item[{\rm (ii)}] If the
origin of system {\rm (\ref{eqnil})} is a focus, then $m \geq 1$,
$m+n$ is even and its cyclicity ${\rm
Cycl}(\mathcal{X}_\varepsilon,p_0)$ satisfies $\, {\rm
Cycl}(\mathcal{X}_\varepsilon,p_0) \, \geq \, (m+n)/2\, -\, 1 $.
In this case, $m$ is the vanishing multiplicity of $V(r,\theta)$
on $r=0$. \item[{\rm (iii)}] If the origin of system {\rm
(\ref{eqnil2})} is a focus and if only analytic perturbations of
$(1,n)$--quasihomogeneous weighted subdegrees $(w_x,w_y)$ with
$w_x \geq n$ and $w_y \geq 2n-1$ are taken into account, then the
maximum number of limit cycles which bifurcate from the origin is
$\lfloor (m-1)/2 \rfloor$, that is, ${\rm
Cycl}(\mathcal{X}_\varepsilon^{[n,2n-1]},p_0)\, = \, \lfloor
(m-1)/2 \rfloor$.
\end{itemize}
\end{theorem}

\begin{remark} \label{remnil} The proof of this theorem shows that if there exists an inverse
integrating factor $V_0^{*}(x,y)$ of system {\rm (\ref{eqnil2})}
such that $V_0^{*}(r \, {\rm Cs}\, \theta, \, r^n \, {\rm Sn}\,
\theta)/r^{2n-1}$ has a Laurent expansion in a neighborhood of
$r=0$, then the exponents of the leading terms of $V_0^{*}(r \,
{\rm Cs}\, \theta, \, r^n  \, {\rm Sn}\, \theta)/r^{2n-1}$ and
$V(r,\theta)$ coincide. Therefore, the value of $m$ can be
determined without performing the transformation of the system to
generalized polar coordinates.
\end{remark}

We assume that the origin of system (\ref{eqnil2}) is a focus with
Andreev number $n$ and that the vanishing multiplicity of an
inverse integrating factor on it is $m$. If system (\ref{eqnil2})
is written as $\dot{x} \, = \, P(x,y)$ and $\dot{y} \, = \,
Q(x,y)$, then the system:
\begin{equation} \label{eqnilp1} \dot{x} \, = \, P(x,y) \, + \, x \, K(x,y,\varepsilon),
\qquad \dot{y} \, = \, Q(x,y) + n y K(x,y,  \varepsilon),
\end{equation} where
\[ K(x,y , \varepsilon)\, = \, \displaystyle \sum_{i=0}^{L-1}
\varepsilon^{L-i} \, a_i \, x^{2i} \] and $L=(m+n)/2-1$, has at
least $(m+n)/2-1$ limit cycles bifurcating from the origin for
convenient values of the real parameters $a_i$. We recall that $m$
and $n$ have the same parity.\newline

The following corollary establishes a necessary condition for
system (\ref{eqnil}) to have an analytic inverse integrating
factor $V_0(x,y)$ defined in a neighborhood of the origin.

\begin{corollary} \label{cornil}
We assume that the origin of system {\rm (\ref{eqnil})} is a
nilpotent focus with Andreev number $n$, and that there exists an
inverse integrating factor $V_0(x,y)$ of {\rm (\ref{eqnil})} which
is analytic in a neighborhood of the origin. Then, $n$ is odd.
\end{corollary}

\subsection{Singular perturbations}

As some recent research papers show, see \cite{DRP}, limit
periodic sets containing an infinite number of critical points may
have a cyclicity higher than expected. Due to the narrow
relationship between limit cycles and the inverse integrating
factor, the context of singular perturbations is a brand new and
very interesting place to apply properties of the inverse
integrating factor in order to detect limit cycles. As we have
seen, for other limit periodic sets, the inverse integrating
factor does not only give an alternative way to study the
cyclicity but contains more information: location of limit cycles,
direct computation of the cyclicity of the object under study, an
explicit partial differential equation (\ref{def-V}) which gathers
all the information, \ldots \newline

As far as we know, the only work where the inverse integrating
factor is related with a singular perturbation problem is
\cite{LlibMedSil}, where one-parameter families of vector fields
$\mathcal{X}_\varepsilon$ in $\mathbb{R}^2$ of the form
$\mathcal{X}_\varepsilon \, = \, f(x,y,\varepsilon)\, \partial_x
\, + \, \varepsilon g(x,y,\varepsilon) \, \partial_y$, where
$\varepsilon \geq 0$ and $f,g$ are analytic functions, are taken
into account. The aim of the singular perturbation problems is to
study the phase portrait, for $\varepsilon$ sufficiently small,
near the set of singular points of $\mathcal{X}_0$, that is,
$\Sigma \, = \, \{(x,y) \in \mathbb{R}^2 \, : \, f(x,y,0)=0\}$. In
particular, the question is to decide if $\mathcal{X}_\varepsilon$
has a limit cycle which tends to a singular orbit of
$\mathcal{X}_0$ when $\varepsilon \searrow 0$. A singular orbit
(also denoted as {\em slow-fast cycle}, see \cite{DRP, DRR}) is a
limit periodic set of the system $\mathcal{X}_0$. For the vector
field $\mathcal{X}_0$, we say that a point $n \in \Sigma$ is {\em
normally hyperbolic} if $(\partial f/\partial x)(n,0) \, \neq \,
0$. The system of differential equations associated to
$\mathcal{X}_\varepsilon$ is
\[ \dot{x} \, = \, f(x,y,\varepsilon), \quad \dot{y} \, = \,
\varepsilon \, g(x,y,\varepsilon),\] where the dot denotes
derivation with respect to the time $t$. We call this system the
{\em fast system}. By the time rescaling $\tau \, = \, \varepsilon
t$, we get the {\em slow system}:
\[ \varepsilon x' \, = \,  f(x,y,\varepsilon), \quad y' \, = \,
g(x,y,\varepsilon), \] where $'$ denotes derivation with respect
to $\tau$. The {\em reduced problem} is defined by the slow system
taking $\varepsilon=0$, which gives one differential equation
constrained to the slow manifold or critical curve $\Sigma$, that
is, the reduced problem is
\[ f(x,y,0)=0, \quad y'\, = \, g(x,y,\varepsilon). \]
The only singular orbits taken into account in \cite{LlibMedSil}
are the ones consisting of three pieces of smooth curves; an orbit
of the reduced problem starting at a normally hyperbolic point
$n_1 \in \Sigma$, an orbit of the reduced problem ending at a
normally hyperbolic point $n_2 \in \Sigma$ and an orbit of the
fast problem connecting the two previous ones. The main results
are the following.
\begin{theorem} {\rm \cite{LlibMedSil}} \label{thLlibMedSil}
Consider $\varepsilon_0>0$ and $V_\varepsilon(x,y)$ an inverse
integrating factor of $X_\varepsilon$, that is
$\mathcal{X}_\varepsilon(V_\varepsilon(x,y)) \, = \, {\rm div}
\mathcal{X}_\varepsilon \, V_\varepsilon(x,y)$, defined in an open
set $\mathcal{U} \subseteq \mathbb{R}^2$ for any $0\leq
\varepsilon \leq \varepsilon_0$. Let $\Gamma \subset \mathcal{U}$
be a singular orbit and $\Gamma_\varepsilon$ be a limit cycle of
$\mathcal{X}_\varepsilon$ in $\mathcal{U}$ for $\varepsilon \in
(0,\varepsilon_0)$, with $\Gamma_\varepsilon \to \Gamma$,
according to the Hausdorff distance. Then $V_0(\Gamma)=0$.
\end{theorem}

\begin{corollary} {\rm \cite{LlibMedSil}} \label{corLlibMedSil}
Consider $V_\varepsilon(x,y)$ an inverse integrating factor of
$X_\varepsilon$ as in Theorem {\rm \ref{thLlibMedSil}}. If the
level zero of the function $V_0(x,y)$ does not contain a closed
curve, then there exists $\varepsilon_0>0$ such that
$\mathcal{X}_\varepsilon$ does not present a limit cycle for
$0<\varepsilon<\varepsilon_0$ in $\mathcal{U}$. \end{corollary}

As an application of these results, the following examples are
given in \cite{LlibMedSil}. The following vector fields present no
limit cycles because the corresponding inverse integrating factors
have no closed curves in their level zero sets.
\begin{itemize}
\item The vector field $\mathcal{X}_\varepsilon \, = \,
(y^2-x^2)\partial_x \, + \,\varepsilon \, x^2\, \partial_y$ has
the inverse integrating factor
$V_\varepsilon(x,y)=y^3-yx^2-x^3\varepsilon$.
\item The vector field $\mathcal{X}_\varepsilon \, = \,
(y-x^2)\partial_x \, + \,\varepsilon \, x\, \partial_y$ has the
inverse integrating factor
$V_\varepsilon(x,y)=-y+x^2+(1/2)\varepsilon$.
\item The vector field $\mathcal{X}_\varepsilon \, = \,
(-y+x^2)\partial_x \, + \,\varepsilon \, x\, \partial_y$ has the
inverse integrating factor
$V_\varepsilon(x,y)=y-x^2+(1/2)\varepsilon$. \end{itemize}

\section{Some generalizations}

\subsection{The inverse Jacobi multiplier}

Inverse Jacobi multipliers are a natural generalization of inverse
integrating factors to $n$-dimensional dynamical systems with $n
\geq 3$. In \cite{BerroneGiacomini2}, it is developed the theory
of inverse Jacobi multiplier from its beginning in the formal
methods of integration of ordinary differential equations to
modern applications.

In this section we will assume that $\mathcal{X} = \sum_{i=1}^n
X_i(x_1,\ldots, x_n) \partial_{x_i}$ is a $\mathcal{C}^1$ vector
field defined in the open set $\mathcal{U} \subseteq
\mathbb{R}^n$. A $\mathcal{C}^1$ function $V : \mathcal{U} \to
\mathbb{R}$ is said to be an inverse Jacobi multiplier for the
vector field $\mathcal{X}$ in $\mathcal{U}$ when $V$ solves in
$\mathcal{U}$ the linear first order partial differential equation
$\mathcal{X} V = V {\rm div} \mathcal{X}$. The first appearance of
these multipliers occurs in the works of C.G.J. Jacobi, about the
middle of the past century. Many properties of inverse integrating factors for the planar case
($n=2$) are inherited by inverse Jacobi multiplier. We list some
of them:
\begin{itemize}
\item If the change of coordinates $y = \phi(x)$ is introduced,
then $W(y) = (V \circ \phi^{-1})(y) \det\{ D \phi(\phi^{-1}(y))
\}$ is an inverse multiplier of the transformed vector field
$\phi_*\mathcal{X}$.

\item Let $V_1$ and $V_2$ be two linearly independent inverse
Jacobi multipliers of $\mathcal{X}$ defined in $\mathcal{U}$. If
$V_1(x) \neq 0$ for all $x \in \mathcal{U}$, then the ratio
$V_2/V_1$ is a first integral of $\mathcal{X}$ in $\mathcal{U}$.

\item One can use local Lie groups of transformations to
find inverse Jacobi multipliers as follows. Assume $\mathcal{X}$
admits in $\mathcal{U}$ a $(n - 1)$--parameter local Lie group of
transformations with infinitesimal generators $\{ \mathcal{Y}_1,
\ldots, \mathcal{Y}_{n-1} \}$. Then, an inverse Jacobi multiplier
$V$ for $\mathcal{X}$ in $\mathcal{U}$ is furnished by the
determinant $V= \det \{ \mathcal{X}, \mathcal{Y}_1, \ldots,
\mathcal{Y}_{n-1} \}$.

\item Let $\{ \mathcal{Y}_1, \ldots, \mathcal{Y}_{n-1} \}$ be
the generators of $n - 1$ local Lie groups of symmetries
admitted by $\mathcal{X}$ in $\mathcal{U}$. Then, the inverse
multiplier $V= \det \{ \mathcal{X}, \mathcal{Y}_1, \ldots,
\mathcal{Y}_{n-1} \}$ vanishes on every invariant solution of
$\mathcal{X}$ contained in $\mathcal{U}$. Recall here that an
invariant solution of $\mathcal{X}$ corresponding to the group $G$
is defined to be an integral curve of $\mathcal{X}$ which is
invariant under the action of $G$.

\item Let $p_0 \in \mathcal{U}$ be an isolated zero of an
inverse Jacobi multiplier $V$ such that $V \geq 0$ in a
neighborhood $\mathcal{N}$ of $p_0$. Then, $p_0$ is a stable
(resp. unstable) singular point of $\mathcal{X}$ provided that
${\rm div} \mathcal{X} \leq 0$ (resp. $\geq 0$) in $\mathcal{N}$.
Furthermore, the stability (resp. unstableness) of $p_0$ is
asymptotic stability (resp. unstableness) provided that ${\rm div}
\mathcal{X} < 0$ (resp. $> 0$) in $\mathcal{N}$.
\end{itemize}

In the following, we summarize some of the results obtained in
\cite{BerroneGiacomini2}. By a {\it limit cycle} $\gamma$ of
$\mathcal{X}$ we mean a $T$--periodic orbit which is $\alpha$ or
$\omega$--limit set of another orbit of $\mathcal{X}$. Let $V$ be
an inverse Jacobi multiplier defined in a region containing
$\gamma$. If $\gamma = \{ \gamma(t) \in \mathcal{U} : 0 \leq t
\leq T \}$, we define
$$
\Delta(\gamma) = \int_0^T {\rm div} \mathcal{X} \circ \gamma(t) \ dt \ .
$$
As it is well known, $\Delta(\gamma)$ is the sum of the
characteristic exponents of the limit cycle $\gamma$. We recall
that if $\Delta(\gamma) > 0$ then $\gamma$ is not orbitally
stable. We will say that $\gamma$ is a {\it strong} limit cycle
when $\Delta(\gamma) \neq 0$. If, on the contrary, $\Delta(\gamma)
= 0$, then we say that $\gamma$ is a {\it weak} limit cycle.

\begin{theorem}{\rm \cite{BerroneGiacomini2}}
Let $V$ be an inverse Jacobi multiplier of $\mathcal{X}$ defined
in a region containing a limit cycle $\gamma$ of $\mathcal{X}$.
Then, $\gamma$ is contained in $V^{-1}(0)$ in the following cases: {\rm (i)} if $\gamma$ is a strong limit cycle, or
{\rm (ii)} if $\gamma$ is asymptotically orbitally stable (unstable).
\end{theorem}

\begin{theorem}{\rm \cite{BerroneGiacomini2}}
Let $V$ be a Jacobi inverse multiplier defined in a neighborhood
of a limit cycle $\gamma$ of the vector field $\mathcal{X}$. If
$\gamma$ is a strong limit cycle, then
\begin{itemize}
\item[{\rm (i)}] $V$ vanishes on $W^s(\gamma)$, the stable manifold of $\gamma$, provided that $\Delta(\gamma) > 0$;
\item[{\rm (ii)}] $V$ vanishes on $W^u(\gamma)$, the unstable manifold of $\gamma$, provided that $\Delta(\gamma) < 0$.
\end{itemize}
\end{theorem}

The following example appears in \cite{BerroneGiacomini2}. Consider
the cubic polynomial vector field in $\mathbb{R}^3$
\begin{equation}
\dot x \, = \, \lambda (-y + x f(x,y)), \quad \dot y \, =\,
\lambda (x + y f(x,y)), \quad \dot z \,=\, z , \label{9-naLJM}
\end{equation}
where $f(x,y) = 1-x^2-y^2$ and $\lambda > 0$ is a real parameter.
The circle $\gamma = \{ f(x,y)= 0 \} \cap \{ z=0 \}$ is a limit
cycle of system (\ref{9-naLJM}) of period $T = 2 \pi / \lambda$.
In fact, $\gamma(t) = (\cos\lambda t, \sin\lambda t, 0)$. It is
easy to compute that
$$
\Delta(\gamma) = \int_0^T {\rm div} \mathcal{X} \circ \gamma(t) dt = \frac{2(1-2 \lambda)}{\lambda} \pi \ .
$$
An inverse Jacobi multiplier for this system is
$$
V_1(x,y,z) = f(x,y) (x^2+y^2) z \ .
$$
In addition, when $\lambda = -1/2$, then $V_2(x,y,z) =
(x^2+y^2)^2$ is another inverse Jacobi last multiplier for system
(\ref{9-naLJM}).
\newline

As usual, a hyperbolic singular point $p_0$ of a $\mathcal{C}^1$
vector field $\mathcal{X}$ is named a {\it saddle point} when the
matrix $D\mathcal{X}(p_0)$ has eigenvalues with both positive and
negative real parts. Assuming that $k$ of these real parts are
positive and the remaining $n - k$ are negative, the stable
manifold theorem ensures the existence of two invariant
$\mathcal{C}^1$ manifolds $W^u(p_0)$ and $W^s(p_0)$ with
dimensions $\dim W^u(p_0) = k$ and $\dim W^s(p_0) = n-k$, such
that they intersect transversally one each other in $p_0$.

\begin{theorem}{\rm \cite{BerroneGiacomini2}}
Let $p_0$ be a nondegenerate strong singular point of the
$\mathcal{C}^1$ vector field $\mathcal{X}$ having an inverse
Jacobi multiplier $V$  defined in a neighborhood of $p_0$. Then
$V$ vanishes on $W^u(p_0)$ {\rm (}resp. $W^s(p_0)${\rm )} provided
that ${\rm div} \mathcal{X}(p_0) < 0$ {\rm (}resp. ${\rm div}
\mathcal{X}(p_0) > 0${\rm )}.
\end{theorem}

\subsection{Time--dependent inverse integrating factors}

In \cite{GaGiMa}, the authors consider autonomous second order differential equations
\begin{equation} \label{Lie-Su-31}
\ddot{x} = w(x, \dot{x}) \ ,
\end{equation}
with $w \in \mathcal{C}^\infty(\mathcal{U})$ and $\mathcal{U}
\subseteq \mathbb{R}^2$ an open set. They associated to
(\ref{Lie-Su-31}) the first order planar system defined on
$\mathcal{U}$ in the usual way
\begin{equation} \label{Lie-Su-32}
\dot{x} = y \ , \ \dot{y} = w(x, y) \ .
\end{equation}
Moreover, it is associated to equations (\ref{Lie-Su-31}) and (\ref{Lie-Su-32})
the vector fields $\mathcal{X} = \partial_t + \dot{x} \partial_x +
w(x, \dot{x}) \partial_{\dot{x}}$ and $\bar{\mathcal{X}} = y
\partial_x + w(x, y) \partial_{y}$, respectively. A $\mathcal{C}^1$ nonconstant function
$I(t,x, y)$ is called an {\it invariant} (or non--autonomous first
integral) of system (\ref{Lie-Su-32}) in $\mathcal{U}$ if it is
constant along the solutions of (\ref{Lie-Su-32}). In other words,
$\mathcal{X} I \equiv 0$ must be satisfied in $\mathcal{U}$. Of
course, we can find at most two functionally independent
invariants of (\ref{Lie-Su-32}). Notice that an invariant provides
information about the asymptotic behavior of the orbits.

A symmetry of (\ref{Lie-Su-31}) is a diffeomorphism $\Phi : (t,x)
\mapsto (\bar{t}, \bar{x})$ that maps the set of solutions of
(\ref{Lie-Su-31}) into itself. Therefore, the symmetry condition
for (\ref{Lie-Su-31}) is just $\bar{x}'' = w(\bar{x},
\bar{x}')$, where the prime denotes the derivative $'= d / d
\bar{t}$. When the symmetry is a 1--parameter Lie group of point
transformations $\Phi_\epsilon$, then $\bar{t} = t + \epsilon \xi(t,x) + O(\epsilon^2)$, $\bar{x} = x
+ \epsilon \eta(t,x) + O(\epsilon^2)$,
for $\epsilon$ close to zero, and the vector field $\mathcal{Y} =
\xi(t,x) \partial_t + \eta(t,x) \partial_x$ is called the {\it
infinitesimal generator} of the 1--parameter Lie group of point
transformations $\Phi_\epsilon$. It is well known that the {\it determining equations} for Lie point symmetries can be
obtained from the linearized condition
\begin{equation} \label{Lie-Su-7}
\mathcal{Y}^{[2]} (\ddot{x} - w(x, \dot{x})) = 0 \ \mbox{when} \
\ddot{x} = w(x, \dot{x}) \ ,
\end{equation}
where $\mathcal{Y}^{[2]} = \mathcal{Y} + \eta^{[1]}(t,x,\dot{x})
\partial_{\dot{x}} + + \eta^{[2]}(t,x,\dot{x}, \ddot{x})
\partial_{\ddot{x}}$ is the so--called {\it second prolongation}
of the infinitesimal generator $\mathcal{Y}$ and $
\eta^{[1]}(t,x,\dot{x}) = D_t \eta - \dot{x} D_t \xi$, $
\eta^{[2]}(t,x,\dot{x},\ddot{x}) = D_t \eta^{[1]} - \ddot{x} D_t
\xi$ where $D_t = \partial_t + \dot{x} \partial_x + \ddot{x}
\partial_{\dot{x}}$ is the operator total derivative with respect to
$t$. Of course, since (\ref{Lie-Su-31}) is autonomous, it
always admits the generator $\mathcal{Y} = \partial_t$ of a Lie
point symmetry. Let $\mathcal{L}_r$ denote the set of all infinitesimal generators
of 1--parameter Lie groups of point symmetries of the differential
equation (\ref{Lie-Su-31}). It is known that $\mathcal{L}_r$ is a
finite dimensional real Lie algebra, where we denote $r = \dim
\mathcal{L}_r$. Moreover, for autonomous second order differential equation we have $r \in \{
1,2,3,8 \}$.

For any $\mathcal{Y}_i = \xi_i(t,x) \partial_t + \eta_i(t,x)
\partial_x \in \mathcal{L}_r$, easily one can check that the Lie bracket
$[\mathcal{X}, \mathcal{Y}_i^{[1]}] = \mu_i(t,x,\dot{x})
\mathcal{X}$ where $\mu_i(t,x, \dot{x}) = \mathcal{X} \xi_i$ and
$\mathcal{Y}_i^{[1]} = \mathcal{Y}_i + \eta_i^{[1]}(t,x,\dot{x})
\partial_{\dot{x}}$ is the first prolongation of $\mathcal{Y}$. If
$r \geq 2$, we define the functions
\begin{equation} \label{Lie-Su-51}
V_{ij}(t,x,\dot{x}) = \det\{\mathcal{X}, \mathcal{Y}_i^{[1]},
\mathcal{Y}_j^{[1]} \} = \left| \begin{array}{ccc} 1 & \dot{x} & w(x,\dot{x}) \\
\xi_i(t,x) & \eta_i(t ,x) & \eta_i^{[1]}(t,x,\dot{x}) \\
\xi_j(t,x) & \eta_j(t ,x) & \eta_j^{[1]}(t,x,\dot{x})
\end{array}\right|
\end{equation}
for $i,j \in \{ 1, \ldots, r\}$ with $1 \leq i < j \leq r$. The aim of the work \cite{GaGiMa} is to generalize the concept of inverse integrating factor $V(x,y)$ of system (\ref{Lie-Su-32}) via the
functions $V_{ij}(t,x,y)$ defined in (\ref{Lie-Su-51}). In fact, in the autonomous particular case $\partial V_{ij} /
\partial t \equiv 0$, we get that $V_{ij}$ is just an inverse integrating factor of (\ref{Lie-Su-32}). On the contrary, when
$\partial V_{ij} / \partial t \not\equiv 0$, in \cite{GaGiMa} it is proved that the zero--sets $V^{-1}(0)$ and $V_{ij}^{-1}(0)$ have similar properties. The next result provides the connection between inverse integrating factors of system (\ref{Lie-Su-32}) and the functions $V_{ij}(t,x,y)$.

\begin{proposition} \label{Lie-Su-53}
Assume that system {\rm (\ref{Lie-Su-32})} possesses an
$r$--dimensional Lie point symmetry algebra with $r \geq 2$ and
define the functions $V_{ij}(t,x,\dot{x})$ as in {\rm
(\ref{Lie-Su-51})}.
\begin{itemize}
\item[{\rm (i)}] $V_{ij}$ satisfies the linear partial differential
equation $\mathcal{X} V_{ij} = V_{ij} \ {\rm div} \mathcal{X}$,
where $\mathcal{X} = \partial_t + \dot{x} \partial_x + w(x,
\dot{x}) \partial_{\dot{x}}$.

\item[{\rm (ii)}] If $r \geq 3$ then, the ratio of any two nonzero
$V_{ij}$ is either a constant or an invariant of {\rm
(\ref{Lie-Su-32})}.

\item[{\rm (iii)}] If $V_{ij} \equiv 0$, then $(\eta_i-y
\xi_i)/(\eta_j-y \xi_j)$ is an invariant of system {\rm
(\ref{Lie-Su-32})}.
\end{itemize}
\end{proposition}

The next theorem is about the invariant curves of
$\bar{\mathcal{X}}$ contained in $V_{ij}^{-1}(0)$ and give us an extension of Theorem 9 in \cite{GLV}
for a case with $\partial V_{ij} / \partial t \not\equiv 0$. We
put special emphasis on periodic orbits of (\ref{Lie-Su-32}) of
any kind (isolated and, therefore, limit cycles or non-isolated
and so belonging to a period annulus). Recall here that a limit
cycle $\gamma := \{ (x(t),y(t)) \in \mathcal{U} :
0 \leq t < T \}$ is {\it hyperbolic} if $\oint_\gamma {\rm
div} \bar{\mathcal{X}}(x(t),y(t)) dt \neq 0$. On the other hand, a
$\mathcal{C}^1$ curve $f(x,y)=0$ defined on $\mathcal{U}$ is invariant for $\bar{\mathcal{X}}$ if $\bar{\mathcal{X}} f = K
f$ for some function $K(x,y)$ called {\it cofactor}.

\begin{theorem}{\rm \cite{GaGiMa}} \label{Lie-Su-36}
Let $\mathcal{U} \subset \mathbb{R}^2$ be an open set and assume
that $\ddot{x} = w(x, \dot{x})$ with $w$ smooth in $\mathcal{U}$
admits an $r$--dimensional Lie point symmetry algebra
$\mathcal{L}_r$ with $r \geq 2$. Consider the functions
$V_{ij}(t,x,\dot{x})$ defined in {\rm (\ref{Lie-Su-51})} for $i,j
\in \{ 1, \ldots, r\}$ with $1 \leq i < j \leq r$. Suppose that
$\gamma=(x(t), y(t)) \subset \mathcal{U}$ is a $T$--periodic orbit
of {\rm (\ref{Lie-Su-32})}. Then the next statements hold:
\begin{itemize}
\item[{\rm (i)}] If $V_{ij}(t,x,\dot{x}) = V(x,\dot{x}) \not\equiv 0$,
with $V \in \mathcal{C}^1(\mathcal{U})$, then $V(x,y)$ is an
inverse integrating factor of system {\rm (\ref{Lie-Su-32})} in
$\mathcal{U}$. In particular, if $\gamma$ is a limit cycle, then
$\gamma \subset \{ V(x,y) = 0 \}$.

\item[{\rm (ii)}] If $V_{ij}(t,x,\dot{x}) = F(t) G(x,\dot{x}) \not\equiv
0$ with non--constants $F$ and $G \in \mathcal{C}^1(\mathcal{U})$,
then $\dot{F} = \alpha F$ with $\alpha \in \mathbb{R} \backslash
\{0\}$ and $G(x,y)=0$ is an invariant curve of system {\rm
(\ref{Lie-Su-32})}. Moreover, we have:
\begin{itemize}
\item[{\rm (ii.1)}] If $\gamma \subset \{ G = 0 \}$ and $G$ is analytic
on $\mathcal{U}$, then $G$ is not square--free, i.e., $G(x,y) =
g^n(x,y) u(x,y)$ with a positive integer $n > 1$ and $g$ and $u$
are analytic functions on $\mathcal{U}$ satisfying $\gamma \subset
\{ g = 0 \}$ and $\gamma \not\subset \{ u = 0 \}$.

\item[{\rm (ii.2)}] If $\gamma \not\subset \{ G = 0 \}$ then $\gamma$ is
hyperbolic and $\alpha T = \oint_\gamma {\rm div}
\bar{\mathcal{X}} (x(t), y(t)) dt$.
\end{itemize}
\end{itemize}
\end{theorem}

An immediate consequence is obtained.

\begin{corollary}{\rm \cite{GaGiMa}} \label{Lie-Su-36-1}
Assume that $\ddot{x} = w(x, \dot{x})$, with $w$ smooth in the
open set $\mathcal{U} \subseteq \mathbb{R}^2$, admits an
$r$--dimensional Lie point symmetry algebra $\mathcal{L}_r$ with
$r \geq 2$. Consider the functions $V_{ij}(t,x,\dot{x})$ for $i,j
\in \{ 1, \ldots, r\}$ with $1 \leq i < j \leq r$. If there is one
$V_{ij}(t,x,y) = F(t) G(x,y) \not\equiv 0$ with non--constants $F$
and $G \in \mathcal{C}^1(\mathcal{U})$, then system {\rm
(\ref{Lie-Su-32})} does not have period annulus in $\mathcal{U}$.
\end{corollary}

In the sequel, we concentrate our attention in the 2--dimensional
case $\mathcal{L}_2$. In \cite{GaGiMa} it is proved that, if
$\partial_t \in \mathcal{L}_2$, then the autonomous or separate
time--variable forms of $V_{ij}(t,x,\dot{x})$ given in Theorem
\ref{Lie-Su-36} are the only possibilities. Moreover, defining the
domain of definition of the infinitesimal generators as the
unbounded open strip $\Xi = \{ (t,x) \in \mathbb{R} \times
\mathbb{X} \} \subset \mathbb{R}^2$, one has the following result.

\begin{theorem}{\rm \cite{GaGiMa}} \label{Lie-Su-42}
Assume that $\ddot{x} = w(x,\dot{x})$ with $w$ smooth in
$\mathcal{U} \subset \mathbb{R}^2$ admits a 2--dimensional Lie
point symmetry algebra $\mathcal{L}_2$ spanned by the
$\mathcal{C}^1(\Xi)$ vector fields $\mathcal{Y}_1 = \partial_t$
and $\mathcal{Y}_2$ such that $[\mathcal{Y}_1, \mathcal{Y}_2]=c_1
\mathcal{Y}_1 + c_2 \mathcal{Y}_2$.
\begin{itemize}
\item[{\rm (i)}] If $c_2=0$ and $\mathcal{Y}_2 \in \mathcal{C}^2(\Xi)$, then
$V_{12}(t,x,\dot{x}) = G(x,\dot{x})$ with $G(x,y) = y^2 [c_1+y
\alpha'(x)-\beta'(x)] + \beta(x) w(x,y)$ an inverse integrating
factor of $\bar{\mathcal{X}}$ in $W = \mathcal{U} \cap
\{\mathbb{X} \times \mathbb{R}\}$ provided that $G \not\equiv 0$.
Moreover, for analytic vector fields $\mathcal{Y}_2$ in $\Xi$,
$\bar{\mathcal{X}}$ has no limit cycles in $W$.

\item[{\rm (ii)}] If $c_2 \neq 0$ then, changing the basis of
$\mathcal{L}_2$ such that $[\bar{\mathcal{Y}}_1,
\bar{\mathcal{Y}}_2] = \bar{\mathcal{Y}}_1$, we have that
$\bar{V}_{12}(t,x,\dot{x}) = \exp(c_2 t) \bar{G}(x,\dot{x})$ with
$\bar{G}(x,\dot{x}) = \dot{x} [c_2 \dot{x} \alpha(x)-c_2 \beta(x)+
\dot{x}^2 \alpha'(x)-\dot{x} \beta'(x)] + \beta(x) w(x,\dot{x})$.
In addition, $\partial w/ \partial x \equiv 0$ or $\beta(x)\equiv
0$. If $\bar{G} \not\equiv 0$ and $\mathcal{U}$ is a simply
connected domain, then $\bar{\mathcal{X}}$ has no periodic orbits
in $\mathcal{U}$ and all the $\alpha$ or $\omega$--limit sets of
$\bar{\mathcal{X}}$ are contained in the invariant curve
$\bar{G}(x,y) = 0$ of $\bar{\mathcal{X}}$.
\end{itemize}
\end{theorem}

As an application of these results to polynomial Li\'enard systems, in \cite{GaGiMa} it is proved the next theorem.

\begin{theorem}{\rm \cite{GaGiMa}} \label{Lie-Su-44}
The polynomial Li\'enard differential equation $\ddot{x} + f(x)
\dot{x} + g(x) = 0$ with $f, g \in \mathbb{R}[x]$ having a
$r$--dimensional Lie point symmetry algebra $\mathcal{L}_r$ with
$r \geq 2$ has no limit cycles in $\mathbb{R}^2$.
\end{theorem}

\section*{Acknowledgements}

We would like to thank Professor H\'ector Giacomini, from
Universit\'e de Tours (France), for his useful comments on this
survey and for encouraging us to study the inverse integrating
factor.

\vspace{0.5cm}

{\bf Addresses and e-mails:} \\
$^{\ (1)}$ Departament de Matem\`atica. Universitat de Lleida.
\\ Avda. Jaume II, 69. 25001 Lleida, SPAIN.
\\ {\rm E--mails:} {\tt garcia@matematica.udl.cat}, {\tt
mtgrau@matematica.udl.cat}

\end{document}